\documentclass[11pt]{article}

\usepackage[a4paper,top=3cm,bottom=3.5cm,left=3cm,right=3cm]{geometry}

\usepackage{amssymb}
\usepackage{amsmath}
\usepackage{amsthm}
\usepackage[english]{babel}
\usepackage[latin1]{inputenc}
\usepackage{latexsym}
\usepackage{mathrsfs}
\usepackage{bbm}
\usepackage{graphicx}
\usepackage[hidelinks]{hyperref}
\usepackage{accents}
\usepackage{cases}
\usepackage{url}
\usepackage{color}

\numberwithin{equation}{section}

\newtheorem{Theorem}{Theorem}[section]
\newtheorem{Proposition}[Theorem]{Proposition}
\newtheorem{Definition}[Theorem]{Definition}
\newtheorem{Lemma}[Theorem]{Lemma}

\theoremstyle{Definition}
\newtheorem{remark}[Theorem]{Remark}

\newcommand{\nn}{\mathbb{N}}
\newcommand{\rr}{\mathbb{R}}

\newcommand{\eee}{\mathbb{E}}

\newcommand{\aaa}{\mathcal{A}}

\newcommand{\eps}{\varepsilon}

\newcommand{\R}{\mathbb{R}}

\newcommand{\E}{\mathcal{E}}
\newcommand{\F}{\mathcal{F}}

\DeclareMathOperator{\esssup}{esssup}
\def\esssup_#1{\underset{#1}{\mathrm{ess\,sup\, }}}
\def\essinf_#1{\underset{#1}{\mathrm{ess\,inf\, }}}
\def\argmax_#1{\underset{#1}{\mathrm{arg\,max\, }}}
\def\argmin_#1{\underset{#1}{\mathrm{arg\,min\, }}}

\def \ep{\hbox{ }\hfill$\Box$}
\def \Sum{\displaystyle\sum}

\def \Int{\displaystyle\int}

\def \Inf{\displaystyle\inf}

\def\b1{\bf 1}

\def \I{\mathbb{I}}

\def \R{\mathbb{R}}
\def \M{\mathbb{M}}

\def \E{\mathbb{E}}
\def \F{\mathbb{F}}

\def \P{\mathbb{P}}

\def \S{\mathbb{S}}

\def \Ac{{\cal A}}

\def \Cc{{\cal C}}
\def \Dc{{\cal D}}

\def \Fc{{\cal F}}
\def \Gc{{\cal G}}

\def \Sc{{\cal S}}

\def \Wc{{\cal W}}

\def \ep{\hbox{ }\hfill$\Box$}
 
\def\reff#1{{\rm(\ref{#1})}}
\def\beqs{\begin{eqnarray*}}
\def\enqs{\end{eqnarray*}}
\def\beq{\begin{eqnarray}}
\def\enq{\end{eqnarray}}

\def \trans{^{\scriptscriptstyle{\intercal}}}

\begin{document}

\title{
Linear-quadratic McKean-Vlasov stochastic control problems with random coefficients on finite and infinite horizon, \\ 
and applications
\thanks{This work is part of the ANR project CAESARS (ANR-15-CE05-0024), and also  
supported by FiME (Finance for Energy Market Research Centre) and the ``Finance et D\'eveloppement Durable - Approches Quantitatives'' EDF - CACIB Chair.}
}

\author{
Matteo BASEI
\footnote{IEOR, University of California, Berkeley,  \sf basei at berkeley.edu}
\qquad\quad
Huy\^en PHAM
\footnote{LPMA, Université Paris Diderot and CREST-ENSAE, \sf pham at math.univ-paris-diderot.fr}
}

\maketitle

\begin{abstract} 
We propose a simple and original approach for solving  linear-quadratic mean-field stochastic control problems. We study both finite-horizon and infinite-horizon pro\-blems, and allow notably  some coefficients to be stochastic. Our  method  is based on a suitable extension  of the martingale formulation for verification theorems in control theory. The optimal control involves the solution to a system of Riccati ordinary differential equations and to a linear mean-field backward stochastic differential equation; existence and uniqueness conditions are provided for such a system. Finally, we illustrate our results through  two applications with explicit solutions:  the first one deals with a portfolio liquidation problem with trade crowding, and the second one considers  an economic model of substitutable production goods. 
\end{abstract}
 
\vspace{2mm}


\vspace{2mm}



\section{Introduction}
\label{Sec:intro}

Given a finite horizon $T>0$,  
(we will later extend the results to the infinite-horizon case), 
we fix a filtered probability space $(\Omega,\Fc,\F,\P)$, where $\F$ $=$ $(\Fc_t)_{0 \leq t \leq T}$  satisfies the usual conditions and is the natural filtration of a standard real Brownian motion $W$ $=$ $(W_t)_{0 \leq t \leq T}$, 
augmented with an independent $\sigma$-algebra $\Gc$.   Let $\rho \geq 0$ be a discount factor, and define the set of admissible (open-loop) controls as
\beqs
\aaa &=&  \left\{ \alpha: \Omega \times [0,T] \to \rr^m \text{ s.t.~$\alpha$ is $\F$-adapted  and}\int_0^T e^{-\rho t} \eee[|\alpha_t|^2] dt < \infty \right\}.
\enqs
Given a square-integrable $\Gc$-measurable random variable $X_0$,  
and a control  $\alpha \in \aaa$, we consider the controlled linear 
mean-field stochastic differential equation in $\rr^d$ defined by
\begin{equation}
\label{pb:SDE}
\begin{cases}
dX^\alpha_t = b_t\big(X^\alpha_t, \eee[X^\alpha_t], \alpha_t, \eee[\alpha_t]\big) dt + \sigma_t\big(X^\alpha_t, \eee[X^\alpha_t], \alpha_t, \eee[\alpha_t]\big) dW_t, \quad 0 \leq t \leq T, 
\\
X^\alpha_0=X_0,
\end{cases}
\end{equation}
where for each $t \in [0,T]$, $x,\bar x \in \rr^d$ and $a,\bar a \in \rr^m$ we have set
\begin{equation}
\label{pb:coeffSDE}
\begin{array}{ccc}
b_t\big(x, \bar x, a, \bar a \big) & = &  \beta_t + B_t x + \tilde B_t \bar x + C_t a + \tilde C_t \bar a, \\
\sigma_t\big(x, \bar x, a, \bar a \big) & = &  \gamma_t + D_t x + \tilde D_t \bar x + F_t a + \tilde F_t \bar a.
\end{array}
\end{equation}
Here, the coefficients $\beta,\gamma$ of the affine terms are vector-valued $\F$-progressively measurable processes, whereas the other coe\-fficients of the linear terms are deterministic matrix-valued processes, see Section \ref{Sec:assumpt} for precise assumptions.  The quadratic cost functional to be minimized over $\alpha$ $\in$ $\Ac$ is 
\begin{equation} \label{pb:payoff}
\begin{array}{ccl}
J(\alpha) & = &  \E \Big[ \Int_0^T e^{-\rho t} f_t\big(X^\alpha_t, \eee[X^\alpha_t], \alpha_t, \eee[\alpha_t]\big) dt + e^{-\rho T}
g\big(X^\alpha_T, \eee[X^\alpha_T] \big) \Big], \\ \\
\rightarrow \;\;\;  V_0 & = &  \Inf_{\alpha \in \aaa} J(\alpha),
\end{array}
\end{equation}
where, for each $t \in [0,T]$, $x,\bar x \in \rr^d$ and $a,\bar a \in \rr^m$ we have set
\begin{equation} \label{pb:coeffPayoff}
\begin{array}{ccl}
f_t\big(x, \bar x, a, \bar a \big) &=& (x-\bar x)\trans Q_t (x-\bar x) + \bar x\trans(Q_t + \tilde Q_t) \bar x +   2a\trans I_t (x-\bar x) + 2 \bar a\trans (I_t+\tilde I_t) \bar x   \\
& & \;\;\; + \;  (a - \bar a)\trans N_t (a - \bar a)  + \bar a\trans(N_t + \tilde N_t) \bar a +  2M_t\trans x + 2H_t\trans a, \\ 
g\big(x, \bar x\big) &=&  (x - \bar x)\trans P (x - \bar x) + \bar x\trans( P+ \tilde P) \bar x + 2L\trans x.  
\end{array}
\end{equation}
Here, the coefficients of the linear terms, $M,H,L$ are vector-valued $\F$-progressively measurable processes, whereas the other coefficients are deterministic matrix-valued processes. 
We refer again  to Section \ref{Sec:assumpt} for the precise assumptions. 
The symbol $\trans$ denotes the transpose of any vector or matrix. 

\vspace{1mm}

\begin{remark}
{\rm {\bf 1.} We have centred in \reff{pb:coeffPayoff}  the quadratic terms in the payoff functions $f$ and $g$. One could equivalently formulate the quadratic terms in non-centred form as    
\begin{equation*} \label{pb:coeffPayoff2}
\begin{array}{ccl}
\tilde f_t\big(x, \bar x, a, \bar a \big) &=& x\trans Q_t x + \bar x\trans \tilde Q_t \bar x  + a\trans N_t a 
+ \bar a\trans \tilde N_t \bar a + 2M_t\trans x  + 2H_t\trans a + 2a\trans I_t x + 2 \bar a\trans \tilde I_t \bar x, \\ 
\tilde g\big(x, \bar x\big) &=&  x\trans P x + \bar x\trans \tilde P \bar x + 2L\trans x,
\end{array}
\end{equation*}
by noting that $\E[ f_t\big(X^\alpha_t, \eee[X^\alpha_t], \alpha_t, \eee[\alpha_t]\big) ]$ $=$ $\E[ \tilde f_t\big(X^\alpha_t, \eee[X^\alpha_t], \alpha_t, \eee[\alpha_t]\big) ]$, 
$\E[g(X^\alpha_T, \E[X^\alpha_T])]$ $=$ $\E[\tilde g(X^\alpha_T, \E[X^\alpha_T])]$, since $Q$, $P$,  $N$, $I$  are assumed to be deterministic.

\vspace{1mm}

\noindent {\bf 2.} Notice that the only coefficients allowed to be stochastic are $\beta,\gamma,M,H,L$. We also note that in \eqref{pb:payoff}-\reff{pb:coeffPayoff}, we do not need to add  a term in the form $\tilde M_t\trans \bar x$: indeed, since we consider the expectation of the running cost, we could equivalently substitute such a term with $\eee[\tilde M_t]\trans x$ by noting that 
$\E[ \tilde M_t\trans \E[X_t^\alpha]]$ $=$ $\E[ \E[\tilde M_t]\trans X_t^\alpha]$. 
Similarly, we do not need to consider  terms $\tilde H_t\trans \bar a$ and $\bar a\trans \check I_t x$, $a\trans \check I_t \bar x$ (for a deterministic matrix $\check I_t$). 
}
\qed
\end{remark}

\vspace{3mm}

Our main result is to prove, by means of a suitable martingale optimality principle,  that there exists, under mild assumption on the coefficients, 
a unique optimal control for the linear-quadratic McKean-Vlasov (LQMKV)  control problem \eqref{pb:payoff}, given by
\beq
\label{introOpt}
\alpha^*_t &=&  - S_t^{-1}U_t(X^*_t - \eee[X^*_t]) -  \hat S_t^{-1} (V_t \eee[X^*_t]+ O_t) - S^{-1}_t (\xi_t-\eee[\xi_t]) , \;\;  t \in [0,T], \hspace{7mm}
\enq
where $X^*_t=X^{\alpha^*}_t$, $S_t,\hat S_t,U_t,V_t,O_t$ are deterministic coefficients, and $(\xi_t)_t$ is a stochastic process defined in \eqref{coeffTHM} of Theorem \ref{thm:optimal}. We see 
notably the effect of the randomness of some coefficients in the model through the last term in \reff{introOpt}.

Next, we extend the results to the infinite-horizon case. In this case all the coefficients are assumed to be constant, except $\beta, \gamma,M, H$, which can be 
time-dependent and stochastic. 
The optimal control is still in the form \eqref{introOpt}, but additional conditions on the coefficients  are required  
in order to have a well-defined problem and an admissible optimal control. 

In recent years, control of McKean-Vlasov equations has gained more and more attention, motivated on one hand from stochastic control for a large population in mean-field interaction, and on the other hand from control problems with cost functional  involving nonlinear functional of the law of the state process (e.g.~the mean-variance portfolio selection problem or risk measure  in finance). 
The literature is now quite important in this topic and we refer to the recent books by Bensoussan, Frehse and Yam \cite{benetal13} and Carmona and Delarue \cite{cardel18}, and the references therein. In this McKean-Vlasov framework, linear-quadratic (LQ) models provides an important class of solvable applications, and has been studied in many papers, including 
\cite{benetal11}, \cite{Yong2013}, \cite{HuangLiYong}, \cite{pha16}, \cite{gra16}, \cite{Sun}, however  mostly for constant or deterministic coefficients.  Let us mention that  [13] deals with stochastic coefficients but considering a priori closed loop strategies in linear form w.r.t the state and its mean.

The contributions of this paper are the following. 
First, we provide a new solving technique for linear-quadratic McKean-Vlasov (LQMKV) control problems. Such problems  are usually tackled by calculus of variations methods via stochastic maximum principle and decoupling techniques. Instead, we give here a different approach, based on an extended version of the standard martingale formulation for verification theorems in control theory, and valid for both finite-horizon and infinite-horizon problems: this method is closely connected to the dynamic programming principle (DPP), which is valid in the LQ framework once taking into account both the state and the mean of the state, hence restoring time consistency of the problem. Notice that \cite{pha16} also uses a DPP approach but in the Wasserstein space of probability measures, and considering a priori closed-loop controls. Our approach is simpler in the sense that it does not rely on the notion of derivative in the Wasserstein space, and considers the larger class of  open-loop controls.  
We are able to obtain analytical solutions via the resolution of a system of two Riccati equations and the solution to a  linear mean-field backward stochastic differential equation. 

As a second contribution of the paper, we extend the results in \cite{Yong2013} by introducing stochastic coefficients in the problem, namely $\beta,\gamma,M,H,L$. To the best of our knowledge, this is the first time that explicit formulas are provided for McKean-Vlasov control pro\-blems with random coefficients in the payoff. The inclusion of randomness in the some coefficients  is an important point, as it leads to a richer class of models, which is useful for many  applications, see e.g.~the investment problem in distributed generation under a random centralised electricity price studied in  \cite{AidBaseiPham}.  We consider two other detailed applications in this paper: our first example studies  the classical portfolio liquidation problem with a general price process subject  to a permanent price impact generated by a large number of market participants; our second example  deals with an infinite horizon  
model of production of an  exhaustible resource with a large number of producers and general random  price process.

The paper is organized as follows. Section \ref{Sec:assumpt}  presents the precise assumptions on the coefficients of the problems and provides a detailed description of the solving technique. In Section \ref{Sec:procedure} we solve, step by step, the control problem. Some remarks on the assumptions, and extensions are collected in Section \ref{Sec:remarks}. In Section \ref{Sec:infinitepb} we adapt the results to the infinite-horizon case. Finally, two  original examples from economics and finance are studied  in Section \ref{Sec:applications}.

\section{Assumptions and verification theorem}
\label{Sec:assumpt}

Throughout the paper, for each $q \in \nn$ we denote by $\mathbb{S}^q$ the set of $q$-dimensional symmetric matrices. Moreover, for each normed space $(\M, |\cdot|)$ we set
\begin{align*}
L^{\infty}([0,T],\M) &= \bigg\{ \phi : [0,T]\to \M \text{ s.t.~$\phi$ is measurable and $\textstyle \sup_{t \in [0,T]} |\phi_t| < \infty$} \bigg\}, 
\\
L^{2}([0,T],\M) &= \bigg\{ \phi : [0,T]\to \M \text{ s.t.~$\phi$ is measurable and $\int_0^T e^{-\rho t} |\phi_t|^2 dt < \infty$} \bigg\},
\\
L^2_{\Fc_T}(\M) &= \{ \phi : \Omega \to \M : \text{$\phi$ is $\Fc_T$-measurable and $\eee[|\phi|^2] < \infty$} \}, 
\\
\Sc^2_{\F}(\Omega \times [0,T],\M) &= 
\bigg\{ \phi : \Omega \times [0,T] \to \M \text{ s.t.~$\phi$ is $\F$-adapted and $\E[\sup_{0\leq t\leq T} |\phi_t|^2] < \infty$} \bigg\},
\\
L^2_{\F}(\Omega \times [0,T],\M) &= \bigg\{ \phi : \Omega \times [0,T]\to \M \text{ s.t.~$\phi$ is $\F$-adapted and $\int_0^T e^{-\rho t} \eee[|\phi_t|^2] dt < \infty$} \bigg\}.
\end{align*}

We ask the following conditions on the coefficients of the problem to hold in the finite-horizon case.

\begin{itemize}
	\item[\textbf{(H1)}] The coefficients in \eqref{pb:coeffSDE} satisfy:
	\begin{itemize}
	 \item[(i)] 	$\beta,\gamma \in L^2_{\F}(\Omega \times [0,T],\rr^d)$,
	 \item[(ii)]  $B,\tilde B, D, \tilde D \in L^\infty([0,T],\rr^{d\times d})$,   $C,\tilde C, F, \tilde F \in L^\infty([0,T],\rr^{d\times m})$.
	 \end{itemize}
 	\item[\textbf{(H2)}] The coefficients in \eqref{pb:coeffPayoff} satisfy:
	\begin{itemize}
	\item[(i)] $Q, \tilde Q \in L^\infty([0,T],\mathbb{S}^d)$, $P, \tilde P \in \mathbb{S}^d$,  $N, \tilde N \in L^{\infty}([0,T],\mathbb{S}^m)$, $I, \tilde I \in L^\infty([0,T],\rr^{m\times d})$,
	\item[(ii)] $M \in L^2_{\F}(\Omega \times [0,T],\rr^d)$, $H \in L^2_{\F}(\Omega \times [0,T],\rr^m)$, $L \in L^2_{\mathcal{F}_T}(\rr^d)$,
	\item[(iii)] there exists $\delta > 0$ such that, for each $t \in [0,T]$,
	\begin{align} \nonumber 
	& N_t \; \geq \;  \delta \, \I_m,
	&& P \; \geq \;  0,
	&& Q_t \!-\! I_t\trans N^{-1}_t I_t \; \geq \;  0,
	\end{align}
	\item[(iv)] there exists $\delta > 0$ such that, for each $t \in [0,T]$,
	\begin{align} \nonumber
	& N_t \!+\! \tilde N_t \; \geq \;  \delta \, \I_m,
	&& P \!+\! \tilde P \; \geq \;  0,
	&& (Q_t \!+\! \tilde Q_t) \!-\! (I_t \!+\! \tilde I_t)\trans (N_t \!+\! \tilde N_t)^{-1} (I_t \!+\! \tilde I_t) \; \geq \;  0.
	\end{align}
	\end{itemize}
\end{itemize}

\vspace{1mm}
\begin{remark}
{\rm The uniform positive definite assumption on $N$ and $N+ \tilde N$ is a standard and natural coercive condition when dealing with linear-quadratic control problems. 
We discuss in Section 4 (see Remark \ref{remH2})  alternative assumptions when $N$ and $\tilde N$ may be degenerate. 
}
\qed
\end{remark}

\vspace{2mm}

By {\bf (H1)} and classical results, e.g.~\cite[Prop. 2.6]{Yong2013}, there exists a unique strong solution $X^\alpha = (X^\alpha_t)_{0 \leq t \leq T}$ to the mean-field SDE \eqref{pb:SDE}, which satisfies the standard estimate 
\begin{equation}
\label{pb:estimX}
\eee\big[ \sup_{t \in [0,T]}  |X^\alpha_t|^2\big] \leq C_\alpha \big (1+\eee[|X_0|^2] \big) < \infty,
\end{equation}
where $C_\alpha$ is a constant which depends on $\alpha$ $\in$ $\Ac$ only via $\int_0^T e^{-\rho t} \eee[|\alpha_t|^2] dt$. Also, by {\bf (H2)} and \eqref{pb:estimX}, the LQMKV control problem  \eqref{pb:payoff} is well-defined, in the sense that 
\begin{equation*}
\text{$J(\alpha) \in \rr$, for each $\alpha \in \aaa$.}
\end{equation*}

To solve the LQMKV control problem, we are going to use a suitable verification theorem. Namely, we consider an extended version of the martingale optimality principle usually cited in stochastic control theory: see Remark \ref{rem:verif} for a discussion.

\begin{Lemma}[\textbf{Finite-horizon verification theorem}]
\label{prop:verif}
Let $\{\Wc_t^\alpha, t \in [0,T], \alpha \in \aaa\}$ be a family of $\F$-adapted processes in the form $\Wc_t^\alpha$ $=$ 
$w_t(X^\alpha_t, \eee[X^\alpha_t])$ for some $\F$-adapted random field $\{w_t(x,\bar x), t\in [0,T], x,\bar x \in \R^d\}$ satisfying
\beq \label{growthw}
w_t(x,\bar x) & \leq & C( \chi_t + |x|^2 + |\bar x|^2), \;\;\; t \in [0,T], \; x,\bar x \in \R^d, 
\enq
for some positive constant $C$, and nonnegative process $\chi$ with $\sup_{t\in [0,T]} \E[|\chi_t|]$ $<$ $\infty$,  and such that
\begin{itemize}
\item[(i)] $w_T(x,\bar x)$ $=$  $g(x,\bar x)$, $x,\bar x$ $\in$ $\R^d$; 
\item[(ii)] the map $t$ $\in$ $[0,T]$ $\longmapsto$ $\E[\Sc_t^\alpha]$, with 
$\Sc_t^\alpha$ $=$ $e^{-\rho t} \Wc_t^\alpha + \int_0^t e^{-\rho s}  f_s\big(X^\alpha_s, \eee[X^\alpha_s], \alpha_s, \eee[\alpha_s]\big) ds$, is nondecreasing for all 
$\alpha \in \aaa$;
\item[(iii)]  the map $t$ $\in$ $[0,T]$ $\longmapsto$ $\E[\Sc_t^{\alpha^*}]$ is constant for some $\alpha^* \in \aaa$.
	\end{itemize}
Then, $\alpha^*$ is an optimal control and $\eee[w_0(X_0,\eee[X_0])]$ is the value of the LQMKV control problem  \eqref{pb:payoff}: 
\begin{equation*}
V_0 \; = \;  \eee[ w_0(X_0,\eee[X_0])] \; = \;  J(\alpha^*).
\end{equation*}
Moreover, any other optimal control satisfies the condition (iii). 
\end{Lemma}
\noindent {\bf Proof.} From the growth condition \reff{growthw} and estimation \reff{pb:estimX}, we see that the function 
\beqs
t \in [0,T] & \longmapsto & \E [ \Sc_t^\alpha]  
\enqs
is well-defined for any $\alpha$ $\in$ $\Ac$. By (i) and (ii), we have for all $\alpha \in \aaa$ 
\beqs
\E[w_0(X_0,\E[X_0])] & = & \E[\Sc_0^\alpha] \\ 
& \leq & \E[\Sc_T^\alpha] \; = \;  \E\big[ e^{-\rho T} g(X_T^\alpha,\E[X_T^\alpha]) +  \int_0^T e^{-\rho t}  f_t\big(X^\alpha_t, \eee[X^\alpha_t], \alpha_t, \eee[\alpha_t]\big) dt \big] \\
& & \hspace{1.13cm}  = \; J(\alpha), 
\enqs
which shows that $\E[w_0(X_0,\E[X_0])]$ $\leq$ $V_0$ $=$ $\inf_{\alpha\in\Ac} J(\alpha)$, since $\alpha$ is arbitrary. Moreover, condition (iii)  with $\alpha^*$ shows that 
$\E[w_0(X_0,\E[X_0])]$ $=$ $J(\alpha^*)$, which proves the optimality of $\alpha^*$ with $J(\alpha^*)$ $=$ $\E[w_0(X_0,\E[X_0])]$.  Finally, suppose that $\tilde\alpha$ $\in$ $\Ac$ is another optimal control. Then
\beqs
\E[\Sc_0^{\tilde\alpha}] \; = \;  \E[w_0(X_0,\E[X_0])]  \; = \; J(\tilde\alpha) &=& \E[\Sc_T^{\tilde\alpha}].
\enqs
Since the map $t$ $\in$ $[0,T]$ $\longmapsto$ $\E[\Sc_t^{\tilde\alpha}]$ is nondecreasing, this shows that this map is actually constant, and concludes the proof. 
\qed
 
\vspace{3mm}

The general procedure to apply such a verification theorem consists of  the following three steps.
\begin{itemize}
\item[-] \emph{Step 1}. We guess a suitable parametric expression for the candidate random field $w_t(x,\bar x)$, 
and set for each $\alpha \in \aaa$ and $t \in [0,T]$,
	\begin{equation}
	\label{defS}
	\Sc^\alpha_t = e^{-\rho t} w_t(X^\alpha_t, \eee[X^\alpha_t]) + \int_0^t e^{-\rho s}  f_s\big(X^\alpha_s, \eee[X^\alpha_s], \alpha_s, \eee[\alpha_s]\big) ds.
	\end{equation}
	\item[-] \emph{Step 2}.  We apply It\^o's formula to $\Sc_t^\alpha$, for $\alpha$ $\in$ $\Ac$, and take the expectation to get
\beqs
d\E[\Sc_t^\alpha] &=& e^{-\rho t} \E[ \Dc_t^\alpha ] dt, 
\enqs
for some $\F$-adapted processes $\Dc^\alpha$ with 
\beqs
 \E[ \Dc_t^\alpha]  &=&  \E \big[ - \rho w_t(X_t^\alpha,\E[X_t^\alpha]) + \frac{d}{dt} \E\big[ w_t(X_t^\alpha,\E[X_t^\alpha])\big] + 
f_t(X_t^\alpha,\E[X_t^\alpha],\alpha_t,\E[\alpha_t]) \big]. \nonumber 
\enqs
We then determine  the coefficients  in the random field $w_t(x,\bar x)$  s.t.~condition (i) in Lemma \ref{prop:verif} (i.e., $w_T(.)$ $=$ $g(.)$) is satisfied, and so as to have
\beqs \label{conDc}
 \; \E[\Dc_t^\alpha] \; \geq \; 0,  \; t \geq 0, \forall \alpha\in \Ac, &\mbox{and}&  \; \E[\Dc_t^{\alpha^*}] \; = \; 0, \; t \geq 0, 
\mbox{ for some } \;  \alpha^* \in \Ac,
\enqs
which ensures that the mean optimality principle conditions (ii) and (iii) are satisfied, and then $\alpha^*$ will be the optimal control. This leads to a system of backward ordinary and stochastic differential equations.
 \item[-] \emph{Step 3}. We study the existence of solutions to the system obtained in Step 2, which will also ensure the square integrability condition of $\alpha^*$ in $\Ac$, hence its optimality.  
\end{itemize}

\vspace{1mm}
\begin{remark} \label{rem:verif}
{\rm
The standard martingale optimality principle used in the verification theorem for stochastic control problems, see e.g.~\cite{elk81}, consists in finding 
a family of processes $\{\Wc_t^\alpha, 0\leq t\leq T, \alpha\in\Ac\}$ s.t. 
\begin{itemize}
\item[(ii')] the process $\Sc_t^\alpha$ $=$ 
$e^{-\rho t} \Wc_t^\alpha + \int_0^t e^{-\rho s}  f_s\big(X^\alpha_s, \eee[X^\alpha_s], \alpha_s, \eee[\alpha_s]\big) ds$, $0\leq t\leq T$, is a submartingale for each $\alpha$ $\in$ 
$\Ac$, 
\item[(iii')]  $\Sc^{\alpha^*}$ is a martingale for some $\alpha^*$ $\in$ $\Ac$, 
\end{itemize}
which obviously implies the weaker conditions (ii) and (iii)  in Proposition \ref{prop:verif}. Practically, the martingale optimality conditions (ii')-(iii')  would reduce via the It\^o decomposition of $\Sc^\alpha$ to the condition that $\mathcal{D}_t^\alpha$ $\geq$ $0$, for each $\alpha \in \aaa$, and $\mathcal{D}_t^{\alpha^*}$ $=$ $0$, $0\leq t\leq T$,  for a suitable control  $\alpha^*$. In the classical framework of stochastic control problem without mean-field dependence, one looks for $\Wc_t^\alpha$ $=$ $w_t(X_t^\alpha)$ for some random field $w_t(x)$ depen\-ding only on the state value, and the martingale optimality principle leads to the classical Hamilton-Jacobi-Bellman (HJB) equation (when all the coefficients are non random) or to a stochastic HJB, see \cite{pen92}, in the general random coefficients case. In our context of McKean-Vlasov control problems, one looks for $\Wc_t^\alpha$ $=$ $w_t(X_t^\alpha,\E[X_t^\alpha])$ depending also on the mean of the state value, and the pathwise condition 
on $\Dc_t^\alpha$ would not allow us  to determine a suitable random field $w_t(x,\bar x)$. Instead, we exploit the weaker condition (ii) formulated as a mean condition on $\E[\Dc_t^\alpha]$, and we shall see in the next section how it leads indeed to a suitable characterization of  $w_t(x,\bar x)$. 
}
\qed
\end{remark}

\section{The optimal control}
\label{Sec:procedure}

In this section, we apply the verification theorem in Lemma \ref{prop:verif} to characterize an optimal control for the problem  \eqref{pb:payoff}. We will follow the procedure outlined at the end of Section \ref{Sec:assumpt}. In the sequel, and for convenience of notations, we set
\begin{equation} \label{nothat}
\begin{aligned}{}
&\hat B_t \; = \; B_t + \tilde B_t, \quad\,\,\, &&\hat C_t \; = \; C_t + \tilde C_t,  \quad\,\,\, &&\hat D_t \; = \; D_t + \tilde D_t, \quad\,\,\, &&\hat F_t \; = \; F_t + \tilde F_t, \\
&\hat I_t \; = \; I_t + \tilde I_t, \quad\,\,\, &&\hat N_t \; = \; N_t + \tilde N_t,  \quad\,\,\, &&\hat Q_t \; = \; Q_t + \tilde Q_t, \quad\,\,\, &&\hat P \; = \; P  + \tilde P. 
\end{aligned}
\end{equation}

\vspace{1mm}
\begin{remark}	\label{rem:notations}
{\rm To simplify the notations, throughout the paper we will often denote expectations by an upper bar and omit the dependence on controls. Hence, for example, we will write $X_t$ for $X^\alpha_t$, $\bar X_t$ for $\eee[X^\alpha_t]$, and $\bar\alpha_t$ for $\E[\alpha_t]$. 
}
\qed
\end{remark}

\noindent {\bf Step 1.}  Given the quadratic structure of the cost functional $f_t$ in \eqref{pb:coeffPayoff}, we infer a candidate for the random field  $\{w_t(x,\bar x), t\in [0,T], x,\bar x \in \R^d\}$ in the form:
\beq \label{wquadra}
w_t(x,\bar x) &=& (x - \bar x)\trans K_t (x - \bar x) + \bar x\trans \Lambda_t \bar x + 2 Y_t\trans x + R_t,
\enq
where $K_t,\Lambda_t,Y_t,R_t$ are  suitable processes to be determined later. 
The centering of the quadratic term is a convenient choice, which provides simpler calculations. Actually, since the quadratic coefficients in the payoff \reff{pb:coeffPayoff} are deterministic symmetric matrices, we  look for deterministic symmetric matrices $K,\Lambda$ as well. Moreover, since in the statement of Lemma \ref{prop:verif} we always consider the expectation of $\Wc^\alpha_t$ $=$ $w_t(X_t^\alpha,\E[X_t^\alpha])$, we can assume, w.l.o.g., that $R$ is deterministic.  Given the randomness of the linear coefficients in \reff{pb:coeffPayoff},  
the process $Y$ is considered in general as an $\F$-adapted process. 
Finally, the terminal condition $w_T(x,\bar x)$ $=$ $g(x,\bar x)$ ((i) in Lemma \ref{prop:verif}) determines the terminal conditions satisfied by $K_T,\Lambda_T,Y_T,R_T$. We then search for processes  $(K,\Lambda,Y,R)$ valued in $\S^d\times \S^d\times \R^d\times \R$ in backward form: 
\begin{equation}\label{KLYfirst}
\left\{
\begin{array}{cclcl}
dK_t &=& \dot K_t dt, & 0 \leq t \leq T, & K_T \; = \; P, \\
d\Lambda_t & = & \dot \Lambda_t dt,  &  0 \leq t \leq T, & \Lambda_T \; = \;  \hat P, \\
dY_t & = & \dot Y_t dt + Z^Y_t dW_t, & 0 \leq t \leq T, & Y_T \; = \; L, \\
dR_t &=& \dot R_t dt, & 0 \leq t \leq T, & R_T \; = \;  0,
\end{array}
\right.
\end{equation}
for some deterministic processes $(\dot K,\dot\Lambda,\dot R)$ valued in $\S^d\times\S^d\times\R$, and $\F$-adapted processes $\dot Y, Z^Y$ valued in $\R^d$. 

\vspace{3mm}

\noindent {\bf Step 2.} For $\alpha \in \aaa$ and $t \in [0,T]$, let $\Sc^\alpha_t$ as in \reff{defS}. We have 
\beqs
d\E[\Sc_t^\alpha] &=& e^{-\rho t} \E[ \Dc_t^\alpha ] dt, 
\enqs
for some $\F$-adapted processes $\Dc^\alpha$ with 
\beqs
\E[ \Dc_t^\alpha] &=&  \E \big[ - \rho w_t(X_t^\alpha,\E[X_t^\alpha]) + \frac{d}{dt} \E\big[ w_t(X_t^\alpha,\E[X_t^\alpha])\big] + 
f_t(X_t^\alpha,\E[X_t^\alpha],\alpha_t,\E[\alpha_t]) \big]. \nonumber 
\enqs
We apply the It\^o's formula to $w_t(X_t^\alpha,\E[X_t^\alpha])$, recalling the quadratic form \reff{wquadra} of $w_t$, the equations in \reff{KLYfirst}, and the dynamics (see equation \reff{pb:SDE})
\beqs
d X_t^\alpha &=& [\bar\beta_t + \hat B_t \bar X_t^\alpha + \hat C_t \bar\alpha_t] dt,  \\
d(X_t^\alpha- \bar X_t^\alpha) & = & \big[ \beta_t - \bar\beta_t + B_t (X_t^\alpha- \bar X_t^\alpha) + C_t (\alpha_t - \bar\alpha_t) \big] dt \\
& &  + \; \big[  \gamma_t + D_t(X_t^\alpha-\bar X_t^\alpha) + \hat D_t \bar X_t^\alpha + F_t(\alpha_t-\bar\alpha_t) + \hat F_t \bar\alpha_t \big] dW_t,
\enqs
where we use the upper bar notation for expectation, see Remark \ref{rem:notations}. Recalling the quadratic form \reff{pb:coeffPayoff} of the running cost $f_t$, we obtain, after careful but straightforward computations, that
\beq 
\E[ \Dc_t^\alpha] 
&=&  \E \Big[ (X_t - \bar X_t)\trans (\dot K_t + \Phi_t ) (X_t - \bar X_t) 
+ \bar X_t\trans\big(\dot\Lambda_t + \Psi_t \big) \bar X_t  \nonumber \\
& & \hspace{2cm}   + \; 2\big( \dot  Y_t + \Delta_t  \big)\trans X_t  + \dot R_t - \rho R_t  + \bar \Gamma_t   +    \chi_t(\alpha)  \Big],    \label{expressD} 
\enq
(we omit the dependence in $\alpha$ of $X$ $=$ $X^\alpha$, $\bar X$ $=$ $\bar X^\alpha$), where
\begin{equation} \label{PhiK}
\left\{
\begin{array}{rcl}
\Phi_t &=& - \rho K_t+ K_tB_t +  B_t\trans K_t + D_t\trans K_t D_t + Q_t \; = \; \Phi_t(K_t), \\
\Psi_t &=& - \rho \Lambda_t +  \Lambda_t\hat B_t + \hat B_t\trans\Lambda_t + \hat D_t\trans K_t \hat D_t + \hat Q_t \; = \; \Psi_t(K_t,\Lambda_t), \\
\Delta_t&=&  - \rho Y_t + B_t\trans Y_t  + \tilde B_t\trans\bar Y_t + K_t (\beta_t - \bar\beta_t) + \Lambda_t \bar\beta_t \\
& & \;\;\;\;\;  + \;  D_t\trans K_t \gamma_t + \tilde D_t\trans K_t \bar\gamma_t + D_t\trans Z_t^Y + \tilde D_t\trans\overline{Z_t^Y} + M_t \; = \; \Delta_t(K_t,\Lambda_t,Y_t,\bar Y_t,Z_t^Y,\overline{Z_t^Y}), \\
\Gamma_t &=&   \gamma_t\trans K_t \gamma_t +  2 \beta_t\trans Y_t + 2 \gamma_t\trans Z_t^Y \; = \; \Gamma_t(K_t,Y_t,Z_t^Y),
\end{array}
\right.
\end{equation}
for $t$ $\in$ $[0,T]$, and
\beq 
\chi_t(\alpha) &=&  (\alpha_t - \bar \alpha_t)\trans S_t (\alpha_t - \bar\alpha_t) + \bar\alpha_t\trans \hat S_t  \bar\alpha_t \nonumber  \\
& & \;\;\; + \;  2 \big(  U_t (X_t-\bar X_t)  + V_t \bar X_t  + O_t + \xi_t - \bar\xi_t  \big)\trans \alpha_t,  \label{defchi}
\enq
with the deterministic coefficients $S_t, \hat S_t, U_t, V_t, O_t$ defined, for $t$ $\in$ $[0,T]$, by
\begin{equation} \label{defSUVO}
\left\{
\begin{array}{rcl}
S_t&=&  N_t + F_t\trans K_t F_t \;=\; S_t(K_t), \\
\hat S_t &=& \hat N_t + \hat F_t\trans K_t \hat F_t \;=\; \hat S_t(K_t), \\
U_t &=& I_t + F_t\trans K_t D_t + C_t\trans K_t \;=\; U_t(K_t), \\
V_t &=& \hat I_t + \hat F_t\trans K_t \hat D_t + \hat C_t\trans \Lambda_t \;=\; V_t(K_t,\Lambda_t), \\
O_t &=& \bar H_t + \hat F_t\trans K_t \bar \gamma_t + \hat C_t\trans \bar Y_t + \hat F_t\trans \overline{Z_t^Y} \;=\; O_t(K_t,\bar Y_t, \overline{Z_t^Y}) ,
\end{array}
\right.
\end{equation}
and the stochastic coefficient $\xi_t$ 
of mean $\bar\xi_t$ 
defined, for $t$ $\in$ $[0,T]$, by
\begin{equation} \label{defxi}
\left\{
\begin{array}{ccl}
\xi_t & = & H_t + F_t\trans K_t \gamma_t  + C_t\trans Y_t + F_t\trans Z_t^Y \;=\; \xi_t(K_t,Y_t,Z_t^Y), \\
\bar\xi_t & = & \bar H_t + F_t\trans K_t \bar \gamma_t  + C_t\trans \bar Y_t + F_t\trans \overline{Z_t^Y} \;=\; \bar\xi_t(K_t,\bar Y_t,\overline{Z_t^Y}).  
\end{array}
\right.
\end{equation}
Here, we have suitably rearranged the terms in \reff{expressD} in order to keep only linear terms in $X$ and $\alpha$, by using the elementary observation that $\E[ \phi_t\trans \bar X_t]$ $=$ $\E[\bar \phi_t \trans X_t]$, and $\E[\psi_t\trans\bar \alpha_t]$ $=$ $\E[\bar \psi_t\trans\alpha_t]$ for any vector-valued random variable $\phi_t$, $\psi_t$ of mean $\bar\phi_t$, $\bar\psi_t$.   
 
Next, the key-point is to complete the square w.r.t.~the control $\alpha$  in the process $\chi_t(\alpha)$ defined in  \reff{defchi}. Assuming for the moment that the symmetric matrices $S_t$ and $\hat S_t$ are positive definite in $\S^m$ (this will follow typically from the nonnegativity of the matrix $K$, as checked in Step 3, and conditions (iii)-(iv) in {\bf (H2)}), it is clear that  one can find  an invertible matrix $\Theta_t$ in $\R^{m\times m}$  (which may be not unique) s.t. 
\beqs
\Theta_t S_t \Theta_t\trans &=& \hat S_t,  
\enqs
for all $t$ $\in$ $[0,T]$, and which is also deterministic like $S_t,\hat S_t$. We can then rewrite the expectation of $\chi_t(\alpha)$ as  
\beqs
\E[\chi_t(\alpha)] &=& \E \Big[ \big( \alpha_t  - \bar\alpha_t + \Theta_t\trans\bar\alpha_t - \eta_t)\trans S_t \big( \alpha_t  - \bar\alpha_t + \Theta_t\trans\bar\alpha_t - \eta_t \big) - \zeta_t \Big],
\enqs
where 
\beqs
\eta_t &=& a_t^0(X_t,\bar X_t) + \Theta_t\trans a_t^1(\bar X_t),
\enqs
with $a_t^0(X_t,\bar X_t)$ a centred random variable, and $a_t^1(\bar X_t)$ a deterministic function
\beqs
a_t^0(x,\bar x) \; = \;   - S_t^{-1} U_t (x-\bar x) -  S_t^{-1} (\xi_t-\bar \xi_t),  & & a_t^1(\bar x) \; = \;   -   \hat S_t^{-1} (V_t \bar x + O_t),  
\enqs
and
\beqs
\zeta_t &=&  (X_t-\bar X_t) \trans\big( U_t\trans S_t^{-1} U_t \big)(X_t-\bar X_t) + \bar X_t\trans \big(V_t\trans \hat S_t^{-1}V_t\big) \bar X_t \\
& & \;\;\; + \;   2\big(U_t\trans S_t^{-1}(\xi_t-\bar \xi_t)  + V_t\trans \hat S_t^{-1} O_t \big)\trans X_t  \\
& & \;\;\; + \; (\xi_t - \bar\xi_t)\trans S_t^{-1} (\xi_t - \bar\xi_t)    +   O_t\trans \hat S_t^{-1} O_t.
\enqs 
We can then  rewrite the expectation in \reff{expressD} as 
\beqs 
 & & \E[\Dc_t^\alpha] \\
&=&  \E \Big[ (X_t - \bar X_t)\trans \big(\dot K_t + \Phi_t^0 \big) (X_t - \bar X_t)  + \;  \bar X_t\trans\big(\dot\Lambda_t + \Psi_t^0 \big) \bar X_t  \nonumber \\
& & \;   + \; 2\big( \dot  Y_t + \Delta_t^0 \big)\trans X_t   + \; \dot R_t - \rho R_t  + \overline{\Gamma_t^0}    \\
& & \; + \;   \big( \alpha_t - a_t^0(X_t,\bar X_t)  - \bar\alpha_t + \Theta_t\trans(\bar\alpha_t - a_t^1(\bar X_t))\big)\trans S_t \big( \alpha_t - a_t^0(X_t,\bar X_t)  - \bar\alpha_t + \Theta_t\trans(\bar\alpha_t - a_t^1(\bar X_t)) \big)  \Big],
\enqs
where we set 
\begin{equation} \label{Phi0}
\left\{
\begin{array}{rcl}
\Phi_t^0 &=& \Phi_t -   U_t\trans S_t^{-1} U_t \; = \; \Phi_t^0(K_t), \\
\Psi_t^0 &=& \Psi_t  -   V_t\trans \hat S_t^{-1}V_t \; = \; \Psi_t^0(K_t,\Lambda_t), \\
\Delta_t^0 &=& \Delta_t -  U_t\trans S_t^{-1}(\xi_t-\bar \xi_t)  - V_t\trans \hat S_t^{-1} O_t \;= \; \Delta_t^0(K_t,\Lambda_t,Y_t,\bar Y_t,Z_t^Y,\overline{Z_t^Y}), \\
\Gamma_t^0 &= & \Gamma_t - (\xi_t - \bar\xi_t)\trans S_t^{-1} (\xi_t - \bar\xi_t)    -   O_t\trans \hat S_t^{-1} O_t \; = \;  
\Gamma_t^0(K_t,Y_t,\bar Y_t,Z_t^Y,\overline{Z_t^Y}),
\end{array}
\right.
\end{equation}
and stress the dependence on $(K,\Lambda,Y,Z^Y)$ in view of \reff{PhiK}, \reff{defSUVO}, \reff{defxi}.   
Therefore, whenever 
\beqs
\dot K_t + \Phi_t^0  \; = \;  0, \qquad \dot\Lambda_t + \Psi_t^0  \;=\;  0, & &  
\dot  Y_t + \Delta_t^0 \; = \;  0, \qquad  \dot R_t - \rho R_t  + \overline{\Gamma_t^0} \; = \;  0
\enqs
holds for all $t$ $\in$ $[0,T]$, we have
\beq
 & & \E[\Dc_t^\alpha] \label{Dcalphafin} \\
&=& \E \Big[  \big( \alpha_t - a_t^0(X_t,\bar X_t)  - \bar\alpha_t + \Theta_t\trans(\bar\alpha_t - a_t^1(\bar X_t))\big)\trans S_t \big( \alpha_t - a_t^0(X_t,\bar X_t)  - \bar\alpha_t + \Theta_t\trans(\bar\alpha_t - a_t^1(\bar X_t)) \big)   \Big], \nonumber
\enq
which is nonnegative for all $0\leq t\leq T$, $\alpha$ $\in$ $\Ac$, i.e.,  the process $\Sc^\alpha$ satisfies the condition (ii) of the verification theorem in Lemma \ref{prop:verif}. We are then led to consider the  following system of backward (ordinary and stochastic) differential equations (ODEs and BSDE):
\begin{equation} \label{sysK}
\left\{
\begin{array}{ccl}
dK_t &=&  -  \Phi_t^0(K_t) dt, \;\;\; 0 \leq t \leq T, \; K_T \; = \; P, \\ 
d\Lambda_t &=& -  \Psi_t^0(K_t,\Lambda_t)   dt,  \;\;\; 0 \leq t \leq T, \; \Lambda_T \; = \; P + \tilde P, \\
dY_t &=& - \Delta_t^0(K_t,\Lambda_t,Y_t,\E[Y_t],Z_t^Y,\E[Z_t^Y]) dt + Z_t^Y dW_t, \;\; 0 \leq t \leq T, \; Y_T \; = \; L,  \\
dR_t &=& \big[ \rho R_t - \E[\Gamma_t^0(K_t,Y_t, \E[Y_t],Z_t^Y,\E[Z_t^Y]) \big] dt, \;\;\; 0 \leq t \leq T, \; R_T \; = \; 0. 
\end{array}
\right.
\end{equation}

\begin{Definition}
A solution to the system \reff{sysK} is a quintuple of processes $(K,\Lambda,Y,Z^Y,R)$ $\in$ 
$L^\infty([0,T],\S^d)\times L^\infty([0,T],\S^d)\times \Sc^2_{\F}(\Omega\times [0,T],\R^d)\times L^2_{\F}(\Omega\times [0,T],\R^d)\times L^\infty([0,T],\R)$ 
s.t.~the $\S^m$-valued  processes $S(K)$, $\hat S(K)$ $\in$ $L^\infty([0,T],\S^m)$, are positive definite a.s., and the following relation 
\begin{equation*}
\left\{
\begin{array}{ccl}
K_t &=& P + \int_t^T \Phi_s^0(K_s) ds, \\
\Lambda_t &=& P + \tilde P + \int_t^T \Psi_s^0(K_s,\Lambda_s) ds, \\
Y_t &=& L + \int_t^T \Delta_s^0(K_s,\Lambda_s,Y_s,\E[Y_s],Z_s^Y,\E[Z_s^Y]) ds +  \int_t^T Z_s^Y dW_s, \\
R_t &=& \int_t^T - \rho R_s + \E[\Gamma_s^0(K_s,Y_s, \E[Y_s],Z_s^Y,\E[Z_s^Y]) \big] ds,
\end{array}
\right.
\end{equation*}
holds for all $t$ $\in$ $[0,T]$. 
\end{Definition}

We shall discuss in the next paragraph (Step 3) the existence of a solution to the system of ODEs-BSDE \reff{sysK}. For the moment, we provide the connection between this system and the solution to the LQMKV control problem. 

\begin{Proposition} \label{profini}
Suppose that $(K,\Lambda,Y,Z^Y,R)$ is a solution to the system of ODEs-BSDE \reff{sysK}. Then, the control process
\beqs
\alpha_t^* &=& a_t^0(X_t^*,\E[X_t^*]) + a_t^1(\E[X_t^*]) \\
&=& - S_t^{-1}(K_t) U_t(K_t) (X_t-\E[X_t^*]) -  S_t^{-1}(K_t) \big(\xi_t(K_t,Y_t,Z_t^Y) - \bar \xi_t(K_t,\E[Y_t],\E[Z_t^Y]\big) \\
& & \;\;\;  -  \;  \hat S_t^{-1}(K_t) \big(V_t(K_t,\Lambda_t) \E[X_t^*] +  O_t(K_t,\E[Y_t], \E[Z_t^Y]) \big),
\enqs
where $X^*$ $=$ $X^{\alpha^*}$ is the state process with the feedback control $a_t^0(X_t^*,\E[X_t^*]) + a_t^1(\E[X_t^*])$,  is the optimal control for the LQMKV problem \reff{pb:payoff}, i.e., 
$V_0$ $=$ $J(\alpha^*)$, and we have
\beqs
V_0 &=& \E\big[(X_0 - \E[X_0])\trans K_0 (X_0 - \E[X_0]) \big]  + \E[X_0]\trans \Lambda_0 \E[X_0] + 2 \E[Y_0\trans X_0]  + R_0. 
\enqs 
\end{Proposition}
\noindent {\bf Proof.} Consider a solution $(K,\Lambda,Y,Z^Y,R)$ to the system \reff{sysK}, and let $w_t$ as of the quadratic form \reff{wquadra}. First, notice that $w$ satisfies the growth condition \reff{growthw} as $K,\Lambda,R$ are bounded and $Y$ satisfies a square-integrability condition  in $L^2_{\F}(\Omega\times [0,T],\R^d)$. The terminal condition $w_T(.)$ $=$ $g$ is also satisfied from the terminal condition of the system \reff{sysK}. Next, for this choice of $(K,\Lambda,Y,Z^Y,R)$, the expectation $\E[\Dc_t^\alpha]$ in \reff{Dcalphafin}  is non\-negative for all $t$ $\in$ $[0,T]$, $\alpha$ $\in$ $\Ac$, which means that the process $\Sc^\alpha$  satisfies the condition (ii) of the verification theorem in Lemma \ref{prop:verif}. Moreover, we see that $\E[\Dc_t^{\alpha}]$ $=$ $0$, $0\leq t\leq T$, for some $\alpha$ $=$ $\alpha^*$ if and only if (recall that $S_t$ is positive definite a.s.)
\beqs 
\alpha_t^* - a_t^0(X_t^*,\E[X_t^*])  - \E[\alpha_t^*] + \Theta_t\trans(\E[\alpha_t^*] - a_t^1(\E[X_t^*])) &=& 0, \;\;\; 0 \leq t \leq T.  
\enqs
Taking expectation in the above relation, and recalling that $\E[a_t^0(X_t^*,\E[X_t^*])]$ $=$ $0$, $\Theta_t$ is invertible, we get $\E[\alpha_t^*]$ $=$ $a_t^1(\E[X_t^*])$, and thus 
\beq \label{expressalpha}
\alpha_t^* &=& a_t^0(X_t^*,\E[X_t^*]) + a_t^1(\E[X_t^*]), \;\;\; 0 \leq t \leq T.  
\enq
Notice that $X^*$ $=$ $X^{\alpha^*}$ is solution to a linear McKean-Vlasov dynamics, and satisfies the square-integrability condition $\E[\sup_{0\leq t\leq T}|X_t^*|^2]$ $<$ $\infty$,  which implies in its turn that $\alpha^*$ satisfies the square-integrability condition  $L^2_{\F}(\Omega\times [0,T],\R^m)$, since $S^{-1}$, $\hat S^{-1}$, $U$, $V$ are bounded, 
and $O$, $\xi$ are square-integrable respectively in  $L^2([0,T],\R^m)$ and $L_{\F}^2(\Omega\times [0,T],\R^m)$. Therefore, $\alpha^*$ $\in$ $\Ac$, and we conclude by the verification 
theorem in Lemma \ref{prop:verif} that it is the unique optimal control.
\qed

\vspace{2mm}

\paragraph{Step 3.} Let us now verify under assumptions {\bf (H1)}-{\bf (H2)} the existence and uniqueness of a  solution to the decoupled system in \eqref{sysK}.

\begin{itemize}
	\item[(i)] We first consider the equation for $K$, which is actually a matrix Riccati equation written as:
\begin{equation} \label{eqK}
\left\{
\begin{array}{rcl}
\frac{d}{dt} K_t + Q_t -\rho K_t + K_tB_t+ B_t\trans K_t + D_t\trans K_t D_t & & \\
- (I_t + F_t\trans K_t D_t + C_t\trans K_t)\trans(N_t + F_t\trans K_t F_t)^{-1}(I_t + F_t\trans K_t D_t + C_t\trans K_t) &=& 0, \; t \in [0,T], \\
K_T &=& P.
\end{array}
\right.
\end{equation}
Multi-dimensional Riccati equations are known to be related to control theory. Namely,  \eqref{eqK} is associated  to the standard linear-quadratic stochastic control problem:
\beqs
v_t(x) &=& \inf_{\alpha\in \Ac} \E \Big[ \int_t^T e^{-\rho t} \Big( (\tilde X_s^{t,x,\alpha})\trans Q_s \tilde X_s^{t,x,\alpha} + 2 \alpha_s\trans I_s \tilde X_s^{t,x,\alpha} 
+ \alpha_s\trans N_s\alpha_s \Big) ds \\
& & \hspace{2cm} + \; e^{-\rho T}  (\tilde X_T^{t,x,\alpha})\trans P  \tilde X_T^{t,x,\alpha}  \Big] ,
\enqs
where $\tilde X^{t,x,\alpha}$ is the controlled linear dynamics solution to
\beqs
d\tilde X_s &=& ( B_s \tilde X_s + C_s \alpha_s) ds + (D_s \tilde X_s + F_s \alpha_s) dW_s, \;\;\; t \leq s \leq T, \; \tilde X_t = x. 
\enqs
 By a standard result in control theory (see \cite[Ch.~6, Thm.~6.1, 7,1, 7.2]{YongZhou}, with a straightforward adaptation of the arguments to include the discount factor), 
 under {\bf (H1)}, {\bf (H2)}(i)-(ii),  there exists a unique solution $K$ $\in$ $L^\infty([0,T],\S^d)$ with  $K_t$ $\geq$ $0$  to \eqref{eqK}, provided that
	\begin{equation}
	\label{ipoK}
	P \; \geq \;  0, \qquad Q_t - I_t\trans N_t^{-1} I_t \; \geq \; 0, \qquad N_t \; \geq \;  \delta \, \I_m, \;\;\; 0 \leq t \leq T, 
	\end{equation}
	for some  $\delta>0$, which is true by  {\bf (H2)}(iii), and in this case, we have $v_t(x)$ $=$ $x\trans K_t x$.  Notice also that $S(K)$ $=$ $N+ F\trans KF$ is positive definite. 
	
	\item[(ii)] Given $K$, we now consider the equation for $\Lambda$. Again, this is a matrix  Riccati equation that we rewrite as
	\begin{equation} \label{eqL}
\left\{
\begin{array}{rcl}
\frac{d}{dt} \Lambda_t + \hat Q^K_t -\rho \Lambda_t + \Lambda_t\hat B_t+ \hat B_t \trans  \Lambda_t & & \\
- \big(\hat I^K_t + \hat C_t\trans\Lambda_t \big)\trans(\hat N^K_t)^{-1} \big(\hat I^K_t + \hat C_t\trans \Lambda_t \big)  &=& 0, \; t \in [0,T], \\
\Lambda_T &=& \hat P ,
\end{array}
\right.
\end{equation}
where we have set, for $t \in [0,T]$,
\beqs
\hat Q_t^K &=& \hat Q_t   + \hat D_t \trans K_t \hat D_t, \\
\hat I^K_t &=&  \hat I_t + \hat F_t \trans K_t \hat D_t, \\
\hat N^K_t &=&   \hat N_t + \hat F_t  \trans K_t \hat F_t.
\enqs
As for the equation for $K$, there exists a unique solution $\Lambda$ $\in$ $L^\infty([0,T],\S^d)$ with  $\Lambda_t$ $\geq$ $0$  to \eqref{eqL}, provided that
	\begin{equation}
	\label{ipoLambda}
	\hat P \; \geq \;  0, \;\;\; \hat Q^K_t - (\hat I^K_t)\trans (\hat N^K_t)^{-1} (\hat I^K_t) \; \geq \;  0, \;\;\;  \hat N^K_t \; \geq \;  \delta \, \I_m, \;\;\; 0 \leq t \leq T, 
	\end{equation}
	for some  $\delta>0$.  
Let us check that {\bf (H2)}(iv) implies \eqref{ipoLambda}. We already have $\hat P$ $\geq$ $0$.  Moreover, as $K\geq 0$ we have: 
$\hat N^K_t$ $\geq$ $\hat N_t$ $\geq$ $\delta \I_m$.  By easy algebraic manipulations  and since $\hat N_t > 0$, we have (omitting the time dependence)
\begin{align*}
	\hat Q^{K} - (\hat I^{K})\trans(\hat N^{K})^{-1} \hat I^{K} & = \hat Q - \hat I\trans \hat N^{-1} \hat I + (\hat D - \hat F \hat N^{-1} \hat I)\trans K (\hat D - \hat F \hat N^{-1} \hat I) 
	\\
	& \qquad - \Big(\hat F\trans K (\hat D - \hat F \hat N^{-1} \hat I) \Big)\trans(\hat N + \hat F\trans K \hat F)^{-1} \Big(\hat F \trans K (\hat D - \hat F \hat N^{-1} \hat I) \Big)
	\\
	& \geq \;  \hat Q - \hat I\trans \hat N^{-1} \hat I + (\hat D - \hat F \hat N^{-1} \hat I)\trans K (\hat D - \hat F \hat N^{-1} \hat I) 
	\\
	& \qquad - \Big(\hat F\trans K (\hat D - \hat F \hat N^{-1} \hat I) \Big)\trans (\hat F\trans K \hat F)^{-1} \Big(\hat F\trans K (\hat D - \hat F \hat N^{-1} \hat I) \Big)
	\\
	& = \; \hat Q - \hat I\trans \hat N^{-1} \hat I  \; \geq \;  0, \;\;\; \mbox{ by {\bf (H2)}(iv)}.
\end{align*} 		
	\item[(iii)] Given $(K,\Lambda)$, we consider the equation for $(Y,Z^Y)$. This is a mean-field linear BSDE written as
 	\begin{equation}
	\label{eqY}
	\begin{cases}
	dY_t = \Big( \vartheta_t + G_t\trans(Y_t - \eee[Y_t]) + \hat G_t\trans \eee[Y_t] + J_t\trans(Z_t^Y - \eee[Z_t^Y]) 
	+ \hat{J}_t\trans\eee[Z_t^Y] \Big) dt + Z_t^Y dW_t, 
	\\
	Y_T = L,
	\end{cases}
	\end{equation}
	where the deterministic coefficients $G$, $\hat G$, $J$, $\hat J$ $\in$ $L^\infty([0,T],\R^{d\times d})$, and the stochastic process  $\vartheta$  $\in$ $L^2_{\F}(\Omega\times [0,T],\R^d)$ 
	are defined by
	\begin{align*}
	& G_t \; = \;  \rho \,\, \I_d  - B_t + C_tS^{-1}_t U_t, 
	\\
	& \hat G_t \; = \;  \rho \,\, \I_d - \hat B_t + \hat C_t\hat S^{-1}_t V_t,
	\\
	& J_t \; = \;  - D_t + F_tS^{-1}_t U_t,
	\\
	& \hat{J}_t \; = \;  - \hat D_t + \hat F_t  \hat S^{-1}_t V_t, 
	\\
	&\vartheta_t \; = \;  - M_t - K_t (\beta_t - \eee[\beta_t] ) - \Lambda_t \eee[\beta_t] - D_t\trans K_t (\gamma_t - \eee[\gamma_t]) - \hat D_t\trans K_t \eee[\gamma_t]
	\\
	& \qquad  \;\;\; + \;  U_t\trans S^{-1}_t\big( H_t - \eee[H_t] + F_t\trans K_t (\gamma_t - \eee[\gamma_t])\big) \;+\;  
	 V_t\trans \hat S^{-1}_t \big(\eee[H_t] + \hat F_t\trans K_t \eee[\gamma_t] \big),
	\end{align*}
and the expressions for $S,\hat S,U,V$ are recalled in \eqref{defSUVO}. By standard results, see \cite[Thm.~2.1]{LiSunXiong}, there exists a unique solution $(Y,Z^Y)$ $\in$ 
$\Sc_{\F}^2(\Omega\times [0,T],\R^d)\times L_{\F}^2(\Omega\times [0,T],\R^d)$ to \eqref{eqY}.  
	\item[(iv)]  Given $(K,\Lambda,Y,Z^Y)$, the equation for $R$ is a linear ODE, whose unique solution is explicitly given  by 
	\beq \label{eqR}
	R_t &=& \int_t^T e^{-\rho (s-t)} h_s ds,
	\enq
	where the deterministic function $h$ is defined, for $t \in [0,T]$, by
	\beqs
	h_t &=&  \eee\big[ - \gamma_t\trans K_t \gamma_t - \beta_t\trans Y_t - 2\gamma_t\trans Z_t^Y  + \xi_t\trans S^{-1}_t \xi_t \big] 
	- \eee[\xi_t]\trans S^{-1}_t \eee[\xi_t] + O_t\trans \hat S^{-1}_t O_t,
	\enqs
	and the expressions of $O$ and $\xi$ are recalled in \reff{defSUVO} and \reff{defxi}. 
\end{itemize}

\vspace{2mm}
To sum up the arguments of this section, we have proved the following result.

\begin{Theorem}
	\label{thm:optimal}
	Under assumptions  {\bf (H1)}-{\bf (H2)}, the optimal control for the LQMKV  pro\-blem  \eqref{pb:payoff} is given  by
	\beq
		\label{optimal}
		\alpha^*_t &=&  - S_t^{-1} U_t(X^*_t - \eee[X^*_t]) - \hat S_t^{-1} (V_t \eee[X^*_t]+ O_t) - S^{-1}_t (\xi_t-\eee[\xi_t]),
	\enq
	where $X^*=X^{\alpha^*}$ and the deterministic coefficients $S,\hat S$  $\in$ $L^\infty([0,T],\S^m)$, $U,V$ in $L^\infty([0,T],\R^{m\times d})$, $O$ $\in$  
	$L^\infty([0,T],\R^{m})$ and the stochastic coefficient $\xi$ $\in$ $L^2_{\F}(\Omega\times [0,T],\R^m)$ are defined by
	\begin{equation}
		\label{coeffTHM}
		\left\{
		\begin{array}{ccl}
			S_t &=& N_t + F_t\trans K_t F_t,
			\\
			\hat S_t &=& N_t + \tilde N_t + (F_t + \tilde F_t)\trans K_t (F_t + \tilde F_t),
			\\
			U_t & =& I_t + F_t\trans K_t D_t + C_t\trans K_t,
			\\
			V_t & =& I_t+\tilde I_t + (F_t + \tilde F_t)\trans K_t (D_t + \tilde D_t) + (C_t+\tilde C_t)\trans \Lambda_t,
			\\
			O_t & =& \eee[H_t] + (F_t + \tilde F_t)\trans K_t \eee[\gamma_t] + (C_t + \tilde C_t)\trans  \eee[Y_t] + (F_t + \tilde F_t)\trans \eee[Z_t^Y],
			\\
			\xi_t & =&  H_t + F_t\trans K_t \gamma_t + C_t\trans Y_t + F_t\trans Z_t^Y,
		\end{array}
		\right. 
	\end{equation}
	with  $(K,\Lambda,Y,Z^Y,R)$ $\in$ 
	$L^\infty([0,T],\S^d)\times L^\infty([0,T],\S^d)\times \Sc^2_{\F}(\Omega\times [0,T],\R^d)\times L^2_{\F}(\Omega\times [0,T],\R^d)\times L^\infty([0,T],\R)$ the unique  solution to \reff{eqK}, \reff{eqL}, \reff{eqY}, \reff{eqR}. The corresponding value of the problem is 
	\beqs
		V_0 &=& J(\alpha^*) \; = \;  \eee\big[ (X_0- \eee[X_0])\trans K_0 (X_0- \eee[X_0]) \big] +  \eee[X_0]\trans \Lambda_0  \eee[X_0] + 2 \eee\big[Y_0\trans X_0] + R_0.
	\enqs
\end{Theorem}

\section{Remarks and extensions}
\label{Sec:remarks}

We collect here some remarks and extensions for the problem presented in the previous sections.

\vspace{1mm}
\begin{remark} \label{remH2}
{\rm  Assumptions {\bf (H2)}(iii)-(iv) are used only for ensuring the existence of a nonnegative solution $(K,\Lambda)$ to equations \reff{eqK}, \reff{eqL}. In some specific cases, they can be substituted by alternative conditions. For example, in the one-dimensional case $n$ $=$ $m$ $=$ $1$ (real-valued control and state variable), with $N$ $=$ $0$ (no quadratic term on the control in the running cost), $I$ $=$ $0$, we see that the equation for $K$ is a first-order linear ODE, which clearly admits a unique solution. Moreover, when $P$ $>$ $0$, then $K$ $>$ $0$ by classical maximum principle, so that, when $F_t$ $\neq$ $0$, we have $S_t$ $>$ $0$. Hence, an alternative condition to {\bf (H2)}(iii) is, for $t \in [0,T]$,

\vspace{4mm}

\noindent {\bf (H2)}(iii') \hspace{1.7cm}  $n$ $=$ $m$ $=$ $1$, $N_t$ $=$ $I_t$ $=$ $0$, $P$ $>$ $0$, $F_t$ $\neq$ $0$. 

\vspace{4mm}

\noindent Let us now discuss an alternative condition to the  uniform positive condition on $N+\tilde N$ in {\bf (H2)}(iv), in the case where $N+\tilde N$ is only assumed to be nonnegative. 
When the constant matrix $P$ is positive definite, then $K$ is uniformly positive definite in $\S^d$, i.e., $K_t$ $\geq$ $\delta \I_m$, $0\leq t\leq T$, for some $\delta$ $>$ $0$, by strong maximum principle for the ODE \reff{eqK}. Then, when $F+\tilde F$ is  uniformly nondegenerate, i.e., $|F_t+\tilde F_t|$ $\geq$ $\delta$, $0\leq t\leq T$, for some $\delta$ $>$ $0$, we see that $\hat S_t$ $=$ $\hat N_t^K$ $\geq$ $(F_t+\tilde F_t)\trans K_t (F_t+\tilde F_t)$ $\geq$ $\delta' \I_d$ for some $\delta'$ $>$ $0$. Notice also that when $I+\tilde I$ $=$ $0$, then $\hat Q^K - (\hat I^K)\trans (\hat N^K)^{-1} (\hat I^K)$ $\geq$ $Q+\tilde Q$.  Consequently, assumption {\bf (H2)}(iv) can be alternatively replaced by 

\vspace{4mm}

\noindent {\bf (H2)}(iv') \hspace{1cm} $N_t+\tilde N_t$,  $P+\tilde P$, $Q_t+\tilde Q_t$ $\geq$ $0$,  $P$ $>$ $0$, $I_t+\tilde I_t$ $=$ $0$,  $|F_t+\tilde F_t|$ $\geq$ $\delta$,  

\vspace{4mm}

\noindent for $t \in [0,T]$ and some $\delta$ $>$ $0$, which ensures that condition \reff{ipoLambda} is satisfied, hence giving  the existence and uniqueness of a nonnegative solution $\Lambda$ to \reff{eqL}. 
}
\qed
\end{remark}

\vspace{1mm}

\begin{remark}  \label{remH22}
{\rm Condition {\bf (H2)} ensures the existence of a nonnegative $\S^d$-valued process $K$ (resp. $\Lambda$) to the Riccati equation \reff{eqK} (resp. \reff{eqL}). In some applications (see for example the ones detailed the next section), this condition is not satisfied, typically as $Q$ $=$ $\tilde Q$ $=$ $0$, while $I$ or $\tilde I$ is nonzero. However, a solution $(K,\Lambda)$ (possibly nonpositive) to these equations may still exist, with $S(K)$ and $\hat S(K)$ positive definite, and one can then still apply Proposition \ref{profini} to get the conclusion of Theorem \ref{thm:optimal}, i.e., the optimal control exists and is given by \reff{optimal}.}
\qed
\end{remark} 

\vspace{1mm}

\begin{remark} \label{remWmulti}
{\rm
The result in Theorem \ref{thm:optimal} can be easily extended to the case where several Brownian motions are present in the controlled equation:
\beqs
dX^\alpha_t &=& b_t\big(X^\alpha_t, \eee[X^\alpha_t], \alpha_t, \eee[\alpha_t]\big) dt 
+ \Sum_{i=1}^n\sigma^i_t\big(X^\alpha_t, \eee[X^\alpha_t], \alpha_t, \eee[\alpha_t]\big) dW^i_t,
\enqs
where $W^1, \dots, W^n$ are standard independent real Brownian motions and, for each $t \in [0,T]$, $i\in\{1,\dots,n\}$, $x,\bar x \in \rr^d$ and $a,\bar a \in \rr^m$, we have set
	\begin{equation}
	\label{newsde}
	\begin{gathered}
	b_t\big(x, \bar x, a, \bar a \big) = \beta_t + B_t x + \tilde B_t \bar x + C_t a + \tilde C_t \bar a,
	\\
	\sigma^i_t\big(x, \bar x, a, \bar a \big) = \gamma^i_t + D^i_t x + \tilde D^i_t \bar x + F^i_t a + \tilde F^i_t \bar a.
	\end{gathered}
	\end{equation}
We ask the coefficients in \eqref{newsde} to satisfy a suitable adaptation of {\bf (H1)}: namely, we substitute $D,\tilde D, F, \tilde F$ with $D^i,\tilde D^i, F^i, \tilde F^i$, for $i \in \{1,\dots,n\}$. The cost functional and {\bf (H2)}  are unchanged.
	
The statement of Theorem \ref{thm:optimal} is easily adapted to this extended framework. To simplify the notations we use Einstein convention: for example, we write 
$(D^i_t)\trans K D^i_t$ instead of $\sum_{i=1}^n (D^i_t)\trans K D^i_t$. The optimal control $\alpha^*$ is given by \eqref{optimal}, where the coefficients are now defined by  
\begin{equation*}
		\left\{
		\begin{array}{ccl}
			S_t &=& N_t + (F_t^i)\trans K_t F_t^i,
			\\
			\hat S_t &=& N_t + \tilde N_t + (F_t^i + \tilde F_t^i)\trans K_t (F_t^i + \tilde F_t^i),
			\\
			U_t & =& I_t + (F_t^i)\trans K_t D_t^i + C_t\trans K_t,
			\\
			V_t & =& I_t+\tilde I_t + (F_t^i + \tilde F_t^i)\trans K_t (D_t^i + \tilde D_t^i) + (C_t+\tilde C_t)\trans \Lambda_t,
			\\
			O_t & =& \eee[H_t] + (F_t^i + \tilde F_t^i)\trans K_t \eee[\gamma_t] + (C_t + \tilde C_t)\trans  \eee[Y_t] + (F_t^i + \tilde F_t^i)\trans \eee[Z_t^Y],
			\\
			\xi_t & =&  H_t + (F_t^i)\trans K_t \gamma_t + C_t\trans Y_t + (F_t^i)\trans Z_t^Y,
		\end{array}
		\right. 
	\end{equation*}
with  $(K,\Lambda,Y,Z^Y,R)$ $\in$ $L^\infty([0,T],\S^d)\times L^\infty([0,T],\S^d)\times \Sc^2_{\F}(\Omega\times [0,T],\R^d)\times L^2_{\F}(\Omega\times [0,T],\R^d)\times L^\infty([0,T],\R)$ the unique  solution to \reff{sysK} with 
\begin{equation} \nonumber
\left\{
\begin{array}{rcl}
\Phi_t(K_t) &=& - \rho K_t + K_tB_t +  B_t\trans K_t + (D_t^i)\trans K_t D_t^i + Q_t ,\\
\Psi_t(K_t,\Lambda_t) &=& - \rho \Lambda_t +  \Lambda_t(B_t+\tilde B_t) + (B_t+\tilde B_t)\trans\Lambda_t,  \\
& & \;\;\; + \;  (D_t^i+\tilde D_t^i)\trans K_t (D_t^i + \tilde D_t^i) + Q_t + \tilde Q_t, \\
\Delta_t(K_t,\Lambda_t,Y_t,\E[Y_t],Z_t^Y,\E[Z_t^Y]) &=&  - \rho Y_t + B_t\trans Y_t  + \tilde B_t\trans\E[Y_t] +  
k (\beta_t - \bar\beta_t) + \Lambda_t \bar\beta_t   + M_t  \\
& & \;   + \;  (D_t^i)\trans K_t \gamma_t + (\tilde D_t^i)\trans K_t \bar\gamma_t + (D_t^i)\trans Z_t^Y + (\tilde D_t^i)\trans\E[Z_t^Y], \\
\Gamma_t(K_t,Y_t,Z_t^Y) &=&   \gamma_t\trans K_t \gamma_t +  2 \beta_t\trans Y_t + 2 \gamma_t\trans Z_t^Y, \;\;\; t \in [0,T]. \hfill\text{\qed}
\end{array}
\right.
\end{equation}
}
\end{remark}

\vspace{1mm}

\begin{remark}
{\rm The optimal control provided by Theorem \ref{thm:optimal} generalizes known results and standard formulas in control theory. 
\begin{itemize}
\item[-] For example, in the case where
	\begin{equation*}
	I_t=\tilde I_t=\beta_t=\gamma_t=M_t=H_t=L_t=0,
	\end{equation*}
then $Y$ $=$ $Z^Y$ $=$ $0$, $R$ $=$ $0$ (correspondingly, we have $O$ $=$ $\xi$ $=$ $0$).  We thus retrieve the formula in \cite[Thm. 4.1]{Yong2013} for the 
optimal control (recalling the notations in \reff{nothat}):
\beqs
\alpha_t^* &=& - \big(N_t + F_t\trans K_t F_t\big)^{-1}(F_t\trans K_t D_t + C_t\trans K_t)(X^*_t - \eee[X^*_t]) \\
& & \; - \; \big(\hat N_t + \hat F_t  \trans K_t \hat F_t \big)^{-1} \big( \hat F_t\trans K_t \hat D_t + \hat C_t\trans \Lambda_t \big) \eee[X^*_t].
\enqs
\item[-] Consider now the case where all the mean-field coefficients are zero, that is 
	\begin{equation*}
	\label{no-mkv}
	\tilde B_t = \tilde C_t = \tilde D_t = \tilde F_t = \tilde Q_t = \tilde N_t = \tilde P_t \equiv 0.
	\end{equation*}
	Assume, in addition, that $\beta_t=\gamma_t=H_t=M_t=0$. In this case, $K_t$ $=$ $\Lambda_t$ satisfy the same Riccati equation, $Y_t=\hat Y_t=R_t=0$, and we have
	\beqs
	S_t &=& \hat S_t \; = \;  N_t + F_t\trans K_t F_t,	\\
	U_t &=&  V_t \; = \;  I_t + F_t\trans K_t D_t + C_t\trans K_t,	\\
	O &=& \xi \; = \;  0,
	\enqs
	which leads to the well-known formula for classical linear-quadratic control problems (see, e.g.~\cite{YongZhou}):
	\begin{equation*}
	\pushQED{\qed} 
	\alpha^*_t \; = \; - (N_t + F_t\trans K_t F_t)^{-1} ( I_t + F_t\trans K_t D_t + C_t\trans K_t)X^*_t, \;\;\; 0 \leq t\leq T.  \qedhere
	\popQED
	\end{equation*}
\end{itemize} 
}
\end{remark}

\vspace{1mm}

\begin{remark}
{\rm The mean  of the optimal state $X^*$ $=$ $X^{\alpha^*}$ can be computed as the solution of a linear ODE. Indeed,by plugging \eqref{optimal} into \eqref{pb:SDE} and taking expectation, we get
\beqs
\frac{d}{dt} \eee[X^*_t] &=&  \big( B_t+\tilde B_t - (C_t + \tilde C_t) \hat S_t^{-1}V_t \big) \eee[X^*_t] + \big( \eee[\beta_t] - (C_t + \tilde C_t) \hat S_t^{-1} O_t \big),
\enqs
which can be solved  explicitly in the one-dimensional case $d$ $=$ $1$, and expressed as an exponential of matrices in the multidimensional case. 
}
\qed
\end{remark}

\section{The infinite-horizon problem}  \label{Sec:infinitepb}

We now consider an infinite-horizon version of the problem in \eqref{pb:payoff} and adapt the results to this  framework. The procedure is similar to the finite-horizon case, but non-trivial technical issues emerge when dealing with the equations for $(K,\Lambda,Y,R)$ and the admissibility of the optimal control.

\vspace{0.2cm}


On a filtered probability space $(\Omega,\Fc,\F,\P)$ as in Section \ref{Sec:intro} with $\F$ $=$ $(\Fc_t)_{t\geq 0}$,  let $\rho > 0$ be a positive discount factor, and define the set of admissible  controls as
\beqs
\aaa &=&  \left\{ \alpha: \Omega \times\R_+ \to \rr^m \text{ s.t.~$\alpha$ is $\F$-adapted and}\int_0^\infty e^{-\rho t} \eee[|\alpha_t|^2] dt < \infty \right\},
\enqs
while the controlled process is defined on $\rr_+$ by
\begin{equation}
\label{pb:SDEINF}
\begin{cases}
dX^\alpha_t = b_t\big(X^\alpha_t, \eee[X^\alpha_t], \alpha_t, \eee[\alpha_t]\big) dt + \sigma_t\big(X^\alpha_t, \eee[X^\alpha_t], \alpha_t, \eee[\alpha_t]\big) dW_t, \quad t \geq 0, 
\\
X^\alpha_0=X_0,
\end{cases}
\end{equation}
where for each $t \geq 0$, $x,\bar x \in \rr^d$ and $a,\bar a \in \rr^m$ we have now set
\begin{equation}
\label{pb:coeffSDEINF}
\begin{aligned}
& b_t\big(x, \bar x, a, \bar a \big) = \beta_t + B x + \tilde B \bar x + C a + \tilde C \bar a,
\\
& \sigma_t\big(x, \bar x, a, \bar a \big) = \gamma_t + D x + \tilde D \bar x + F a + \tilde F \bar a.
\end{aligned}
\end{equation}
Notice that, unlike Section \ref{Sec:intro} and as usual in infinite-horizon problems, the coefficients of the linear terms are constant vectors, but  the coefficients 
$\beta$ and $\gamma$ are allowed to be stochastic processes. 

\vspace{1mm}

The control problem on infinite horizon is formulated as
\begin{equation} \label{pb:payoffINFI}
\begin{array}{ccl}
J(\alpha) &=&  \eee \Big[ \int_0^\infty e^{-\rho t}  f_t\big(X^\alpha_t, \eee[X^\alpha_t], \alpha_t, \eee[\alpha_t]\big) dt \Big], \\
 V_0 &=&  \Inf_{\alpha \in \aaa}  J(\alpha),
\end{array}
\end{equation}
where, for each $t \geq 0$, $x,\bar x \in \rr^d$ and $a,\bar a \in \rr^m$ we have set
\beq 
f_t\big(x, \bar x, a, \bar a \big) &=& (x-\bar x)\trans Q (x-\bar x) + \bar x\trans(Q + \tilde Q) \bar x +   2a\trans I (x-\bar x) + 2 \bar a\trans (I+\tilde I) \bar x   \nonumber  \\
& & \;\;\; + \;  (a - \bar a)\trans N (a - \bar a)  + \bar a\trans(N + \tilde N) \bar a +  2M_t\trans x + 2H_t\trans a. \label{pb:coeffPayoffINF}
\enq

Notice that, as usual in  infinite-horizon problems, the coefficients of the quadratic terms are here constant matrices, and the only non-constant coefficients are $H,M$, which may be stochastic processes. Given a normed space $(\M,|.|)$, we set 
\beqs
L^{\infty}(\R_+,\M) &=& \bigg\{ \phi : \R_+\to \M \text{ s.t.~$\phi$ is measurable and $\textstyle \sup_{t \geq 0} |\phi_t| < \infty$} \bigg\}, \\
L^{2}(\R_+,\M) &=& \bigg\{ \phi : \R_+ \to \M \text{ s.t.~$\phi$ is measurable and $\textstyle \int_0^\infty e^{-\rho t} |\phi_t|^2 dt  < \infty$} \bigg\}, \\
L^2_{\F}(\Omega \times \R_+,\M) &=& \bigg\{ \phi : \Omega \times \R_+\to \M \text{ s.t.~$\phi$ is $\F$-adapted and $\int_0^\infty e^{-\rho t} \eee[|\phi_t|^2] dt < \infty$} \bigg\},
\enqs
and ask the following conditions on the coefficients of the problem to hold in the infinite-horizon case.

\begin{itemize}
\item[\textbf{(H1')}] The coefficients in \eqref{pb:coeffSDEINF} satisfy:
	\begin{itemize}
	 \item[(i)] 	$\beta,\gamma \in L^2_{\F}(\Omega \times\R_+,\rr^d)$,
	 \item[(ii)]  $B,\tilde B, D, \tilde D \in \rr^{d\times d}$, $C,\tilde C, F, \tilde F \in  \rr^{d\times m}$.
	 \end{itemize}
 \item[\textbf{(H2')}] The coefficients in \eqref{pb:coeffPayoffINF} satisfy:
	\begin{itemize}
	\item[(i)] $Q, \tilde Q \in  \mathbb{S}^d$,   $N, \tilde N \in  \mathbb{S}^m$, $I, \tilde I \in   \rr^{m\times d}$,
	\item[(ii)] $M \in L^2_{\F}(\Omega \times\R_+,\rr^d)$, $H \in L^2_{\F}(\Omega \times\R_+,\rr^m)$,  
	\item[(iii)]  $ N \; > \;  0$, \, $Q \!-\! I\trans N^{-1} I \; \geq \;  0$,
	\item[(iv)]  $ N \!+\! \tilde N \; > \;  0$, \, $(Q \!+\! \tilde Q) \!-\! (I \!+\! \tilde I)\trans (N \!+\! \tilde N)^{-1} (I \!+\! \tilde I) \; \geq \;  0$.
	\end{itemize}
\item[\textbf{(H3')}]  The coefficients in \eqref{pb:coeffSDEINF} satisfy $\rho$ $>$  $2(|B| + |\tilde B| + |D|^2 + |\tilde D|^2)$.  
\end{itemize}

\vspace{2mm}

Assumptions {\bf (H1')} and {\bf (H2')} are simply a rewriting of {\bf (H1)} and {\bf (H2)} for the case where the coefficients do not depend on the time. A further condition {\bf (H3')}, not present in the finite-horizon case, is here required in order to have a well-defined problem, as justified below.

By {\bf (H1')} and classical results, there exists a unique strong solution $X^\alpha = (X^\alpha_t)_{t \geq 0}$ to the SDE in \eqref{pb:SDEINF}. Moreover, by {\bf (H1')} and 
{\bf (H3')}, standard estimates (see Lemma \ref{lem:est} below) leads to the following estimate:
\begin{equation}
\label{pb:estimXINF}
\int_0^\infty e^{-\rho t} \eee[|X^\alpha_t|^2] dt \leq \tilde C_\alpha(1+\eee[|X_0|^2]) < \infty,
\end{equation}
where $\tilde C_\alpha$ is a constant which depends on $\alpha$ $\in$ $\Ac$ only via $\int_0^\infty e^{-\rho t} \eee[|\alpha_t|^2] dt$. Also, by {\bf (H2')} and \eqref{pb:estimXINF}, the problem in \eqref{pb:payoffINFI} is well-defined, in the sense that $J(\alpha)$ is finite for each $\alpha \in \aaa$. 
Even if the proof of \eqref{pb:estimXINF} is an easy adaptation of standard estimates, we prefer to briefly report the arguments here, in order to outline where the expression in 
{\bf (H3')} comes from.

\begin{Lemma}
\label{lem:est}
Under {\bf (H1')} and {\bf (H3')}, the estimate in \eqref{pb:estimXINF} holds for each square-integrable variable $X_0$ and $\alpha \in  \aaa$. 
\end{Lemma}
\noindent {\bf Proof.}
By  It\^o's formula and the definition of $b_t,\sigma_t$ in \reff{pb:coeffSDEINF}, we have (shortened notations, see Remark \ref{rem:notations}) 
\beq
\frac{d}{dt}\eee\big[e^{-\rho t}|X_t|^2\big] &=&  e^{-\rho t}\eee\big[-\rho |X_t|^2 + 2  b_t\trans X_t + |\sigma_t|^2 \big] \nonumber \\
&\leq& e^{-\rho t}\eee\Big[-\rho |X_t|^2 + 2 \Big( \beta_t\trans X_t + X_t\trans BX_t + \bar X_t\trans\tilde B \bar X_t + \alpha_t\trans CX_t + \bar \alpha_t\trans \tilde C X_t \Big)  \nonumber \\
& & \qquad + \;  2 \Big( |\gamma_t|^2+ |D|^2|X_t|^2 + |\tilde D|^2 |\bar X_t|^2 + |F|^2|\alpha_t|^2 + |\tilde F|^2|\bar \alpha_t|^2 \Big) \Big].  \label{estimXrho}
\enq
Then, the Young inequality leads, for each $\eps>0$, to
\begin{align*}
\frac{d}{dt}\eee\big[e^{-\rho t}|X_t|^2\big] \, &\leq e^{-\rho t}\eee\Big[-\rho |X_t|^2 + \eps^{-1}|\beta_t|^2 + \eps|X_t|^2 + 2 |B||X_t|^2 + 2|\tilde B||\bar X_t|^2 \\
&\qquad + \eps^{-1}|C|^2|\alpha_t|^2 + \eps|X_t|^2 + \eps^{-1}|\tilde C|^2|\bar \alpha_t|^2 + \eps|X_t|^2 \\
&\qquad + 2 \Big( |\gamma|^2+ |D|^2|X_t|^2 + |\tilde D|^2 |\bar X_t|^2 + |F|^2|\alpha_t|^2 + |\tilde F|^2|\bar \alpha_t|^2 \Big) \Big]
\\ 
&\leq e^{-\rho t} \Big\{ \big(-\rho + 2|B| + 2|\tilde B| + 2|D|^2 + 2|\tilde D|^2 + 3\eps \big) \eee[|X_t|^2] \\
&\qquad + \big(\eps^{-1}(|C|^2+|\tilde C|^2) + 2(|F|^2 +|\tilde F|^2) \big)\eee[|\alpha_t|^2] + \big(\eps^{-1}\eee[|\beta_t|^2]+2\eee[|\gamma_t|^2]\big) \Big\}.
\end{align*}
We then have, for each $t \geq 0$,
\begin{equation*}
\eee[e^{-\rho t}|X^\alpha_t|^2] \leq \xi^\alpha + \eta \int_0^t \eee[e^{-\rho s}|X^\alpha_s|^2] ds,
\end{equation*}
where we have set
\begin{align*}
& \xi^\alpha = \eee[|X_0|^2] + \int_0^\infty e^{-\rho s} \big(\eps^{-1} \eee[|\beta_s|^2]+2 \eee[|\gamma_s|^2]\big) ds \\
& \qquad\qquad\qquad + \Big(\eps^{-1}(|C|^2 + |\tilde C|^2) + 2 (|F|^2 + |\tilde F|^2) \Big) \int_0^\infty e^{-\rho s} \eee[|\alpha_s|^2] ds,
\\
& \eta = -\rho + 2(|B| + |\tilde B| + |D|^2 + |\tilde D|^2 \big) + 3\eps.
\end{align*}
By Gronwall's lemma, we deduce that
\begin{equation*}
\eee[e^{-\rho t}|X^\alpha_t|^2] \leq \xi^\alpha e^{\eta t},
\end{equation*}	
which concludes the proof as $\eta<0$ for $\eps$ small enough, by {\bf (H3')}.
\qed


\vspace{3mm}

\begin{remark}
\label{RemSimplH3}
{\rm 
In some particular cases, the condition in {\bf (H3')} can be weakened. For example, assume that $B \leq 0$ (similar argument if $\tilde B \leq 0$). Then, since $x\trans Bx \leq 0$ for each $x$, one term in the  above estimate  \reff{estimXrho} can be simplified, so that {\bf (H3')} simply writes
\beqs
\rho &>&  2(|\tilde B| + |D|^2 + |\tilde D|^2).
\enqs
As a second example, assume that $\gamma = F=\tilde F=0$. If, in addition, $\tilde D =0$ (resp.~$D=0$), the condition in {\bf (H3')} simplifies to
	\beqs
	\rho &>&  2|B| + |D|^2 \qquad\quad \text{(resp.~$\rho$ $>$ $2|B| + |\tilde D|^2$)}.
	\enqs
Indeed, since we have $|\sigma|^2=|D|^2|X|^2$ (resp.~$|\sigma|^2=|\tilde D|^2|X|^2$), we do not need the estimate on the square which introduces the factor $2$ in \reff{estimXrho}.
}
\qed
\end{remark}

The infinite-horizon version of the verification theorem is an easy adaptation of the arguments in Lemma \ref{prop:verif}.

\begin{Lemma}[\textbf{Infinite-horizon verification theorem}]
\label{prop:verifINF}
Let $\{\Wc_t^\alpha, t \geq 0, \alpha \in \aaa\}$ be a family of $\F$-adapted processes in the form $\Wc_t^\alpha$ $=$ 
$w_t(X^\alpha_t, \eee[X^\alpha_t])$ for some $\F$-adapted random field $\{w_t(x,\bar x), t\geq 0, x,\bar x \in \R^d\}$ satisfying
\beq \label{growthwinfi}
w_t(x,\bar x) & \leq & C( \chi_t + |x|^2 + |\bar x|^2), \;\;\; t \in \R_+, \; x,\bar x \in \R^d, 
\enq
for some positive constant $C$, and nonnegative process $\chi$ s.t.~$e^{-\rho t} \E[\chi_t]$ converges to zero as $t \to \infty$,  and such that
\begin{itemize}
\item[(i)] the map $t$ $\in$ $\R_+$ $\longmapsto$ $\E[\Sc_t^\alpha]$, with 
$\Sc_t^\alpha$ $=$ $e^{-\rho t} \Wc_t^\alpha + \int_0^t e^{-\rho s}  f_s\big(X^\alpha_s, \eee[X^\alpha_s], \alpha_s, \eee[\alpha_s]\big) ds$, is nondecreasing for all 
$\alpha \in \aaa$;
\item[(ii)]  the map $t$ $\in$ $\R_+$ $\longmapsto$ $\E[\Sc_t^{\alpha^*}]$ is constant for some $\alpha^* \in \aaa$. 
	\end{itemize}
Then, $\alpha^*$ is an optimal control and $\eee[w_0(X_0,\eee[X_0])]$ is the value of the LQMKV control problem  \eqref{pb:payoffINFI}: 
\begin{equation*}
V_0 \; = \;  \eee[ w_0(X_0,\eee[X_0])] \; = \;  J(\alpha^*).
\end{equation*}
Moreover, any other optimal control satisfies the condition (iii). 	
\end{Lemma}
\noindent {\bf Proof.}
Since the integral in \eqref{pb:estimXINF} is finite,  we have $\lim_{t\rightarrow\infty} e^{-\rho t} \eee[|X_t^\alpha|^2]=0$ for each $\alpha$; by \reff{growthwinfi}, we deduce that 
$\lim_{t\rightarrow\infty} e^{-\rho t} \eee[|\Wc^\alpha_t|]$ $=$  $\lim_{t\rightarrow\infty} e^{-\rho t} \eee[|w_t(X^\alpha_t, \eee[X^\alpha_t])|]$ $=$  $0$. Then, the rest of the proof follows the same arguments as in the one of Lemma \ref{prop:verif}.
\qed

\vspace{2mm}

The steps  to apply Lemma  \ref{prop:verifINF} are the same as the ones in the finite-horizon case. We report the main changes in the arguments with respect to the 
finite-horizon case. 

\vspace{2mm}

\noindent {\bf Steps 1-2.}  
We search for a random field $w_t(x,\bar x)$ in a quadratic form as in \reff{wquadra}: 
\beqs
w_t(x,\bar x) &=& (x - \bar x)\trans K_t (x - \bar x) + \bar x\trans \Lambda_t \bar x + 2 Y_t\trans x + R_t,
\enqs
where the mean optimality principle of Lemma \ref{prop:verifINF} leads now to the following system 
\begin{equation} \label{sysKinfi}
\left\{
\begin{array}{ccl}
dK_t &=&  -  \Phi_t^0(K_t) dt, \;\;\; t \geq 0,   \\ 
d\Lambda_t &=& -  \Psi_t^0(K_t,\Lambda_t)   dt,  \;\;\;  t \geq 0,   \\
dY_t &=& - \Delta_t^0(K_t,\Lambda_t,Y_t,\E[Y_t],Z_t^Y,\E[Z_t^Y]) dt + Z_t^Y dW_t, \;\;  t \geq 0,    \\
dR_t &=& \big[ \rho R_t - \E[\Gamma_t^0(K_t,Y_t, \E[Y_t],Z_t^Y,\E[Z_t^Y]) \big] dt, \;\;\; t \geq 0,  
\end{array}
\right.
\end{equation}
Notice that there are no terminal conditions in the system, since we are considering an infinite-horizon framework. The maps $\Phi^0$, $\Psi^0$, $\Delta^0$, $\Gamma^0$ are defined as in \reff{Phi0}, \reff{PhiK}, \reff{defSUVO}, \reff{defxi}, where the coefficients $B,\tilde B, C,\tilde C, D,\tilde D, F,\tilde F, Q, \tilde Q, N, \tilde N, I, \tilde I$ are now constant. 

\vspace{2mm}

\noindent {\bf Step 3.} We now prove the existence of a solution to the system in \reff{sysKinfi}.

\begin{itemize}
\item[(i)] Consider the ODE for $K$. Notice that the map $k$ $\in$ $\S^d$ $\mapsto$ $\Phi^0(k)$ does not depend on time as the coefficients are constant. We then look for a constant 
nonnegative matrix $K$ $\in$ $\S^d$ satisfying $\Phi^0(K)$ $=$ $0$, i.e., solution to
 \beq
 Q -\rho K + KB+ B\trans K + D\trans K D & &  \nonumber \\
-  \;  (I + F\trans K D + C\trans K)\trans(N + F\trans K F)^{-1}(I + F\trans K D + C\trans K) &=& 0.  \label{eqKINF}
 \enq
As in the finite-horizon case, we prove the existence of a solution to \reff{eqKINF} by relating it to a suitable  infinite-horizon linear-quadratic control problem. However, as we could not find a direct result in the literature for such a connection, we proceed through  an approximation argument. For $T \in \rr_+ \cup \{\infty\}$ and $x \in \rr^d$, we consider the following control problems:
\beqs
V^T(x) & : =&  \inf_{\alpha \in \mathcal{A}_T} \eee \Big[ \int_0^T e^{-\rho t} \Big( (\tilde X_t^{x,\alpha})\trans Q \tilde X_t^{x,\alpha} + 2 \alpha_t\trans  I \tilde X_t^{x,\alpha}  
			+ \alpha_t\trans N \alpha_t  \Big) dt \Big],
\enqs
where $\mathcal{A}_T$ $=$ $\Big\{ \alpha: \Omega \times [0,T[ \to \rr^d \text{ s.t.~$\alpha$ is $\F$-adapted  and} \int_0^T e^{-\rho t} \eee[|\alpha_t|^2] dt < \infty \Big\}$, and for $\alpha$ $\in$ $\Ac_T$, 
$(\tilde X^{x,\alpha})_{0\leq t\leq T}$ is the solution to 
\beqs
d\tilde X_t &=& (B \tilde X_t + C \alpha_t) dt + (D \tilde X_t + F \alpha_t) dW_t, \qquad \tilde X_0 \; = \; x.
\enqs
The integrability condition $\alpha \in \mathcal{A}_T$ implies that $\int_0^T e^{-\rho t} \eee[|\tilde X^{x,\alpha}_t|^2] dt < \infty$, and so the problems $V_T(x)$ are well-defined for any 
$T \in \rr_+ \cup \{\infty\}$. If $T$ $<$ $\infty$,  as already recalled in the finite-horizon case, {\bf (H1')}-{\bf (H2')} imply that there exists a unique symmetric 
solution $\{K^T_{t}\}_{t \in [0,T]}$ to
\begin{equation} \label{eqKinter}
\left\{
\begin{array}{rcl}
\frac{d}{dt} K_t^T + Q -\rho K_t^T + K_t^TB+ B\trans K_t^T + D\trans K_t^T D & & \\
- (I + F\trans K_t^T D + C\trans K_t^T)\trans(N + F\trans K_t^T F)^{-1}(I + F\trans K_t^T D + C\trans K_t^T) &=& 0, \; t \in [0,T], \\
K_T^T &=& 0, 
\end{array}
\right.
\end{equation}
and that for every $x \in \rr^d$ we have: $V^T(x)$ $=$ $x\trans K_0^T x$. It is easy to check from the definition of $V^T$ that $V^T(x)$ $\rightarrow$ $V^\infty(x)$ as $T$ goes to infinity, from which we deduce that  
\begin{equation*}
V^\infty(x) \; = \;  \lim_{T\rightarrow\infty}  x\trans K_{0}^T x  \; = \;  x\trans(\lim_{T\rightarrow\infty}  K_{0}^T )x, \;\;\; \forall x \in \R^d. 
\end{equation*}
This implies the existence of the limit $K$ $=$ $\lim_{T\rightarrow\infty} K_0^T$.  By passing to the limit in $T$ in the above ODE \reff{eqKinter} at $t$ $=$ $0$, we obtain by standard arguments (see, e.g., Lemma 2.8 in \cite{SunYong}) that  $K$ satisfies \eqref{eqKINF}.  Moreover, $K \in \S^d$ and $K \geq 0$ as it is the limit of symmetric nonnegative matrices.
\item[(ii)] Given $K$, we now consider the equation for $\Lambda$.  Notice that the map $\ell$ $\in$ $\S^d$ $\mapsto$ $\Psi^0(K,\ell)$ does not depend on time as the coefficients are constant. We then look for a constant nonnegative matrix $\Lambda$ $\in$ $\S^d$ satisfying $\Phi^0(K,\Lambda)$ $=$ $0$, i.e., solution to
\beq
\hat Q^K  -\rho \Lambda  + \Lambda(B+\tilde B) + (B+\tilde B) \trans  \Lambda & & \nonumber \\
- \big(\hat I^K + (C+\tilde C)\trans\Lambda \big)\trans(\hat N^K)^{-1} \big(\hat I^K + (C+\tilde C)\trans \Lambda \big)  &=& 0, \label{eqLINF}
\enq
where we set
\beqs
\hat Q^K &=& (Q + \tilde Q) + (D+\tilde D)\trans K (D+\tilde D), \\
\hat I^K &=& (I+\tilde I) + (F + \tilde F)\trans K (D + \tilde D), \\
\hat N^K &=& (N + \tilde N) + (F + \tilde F)\trans K (F + \tilde F). 
\enqs
Existence of a solution to \eqref{eqLINF} is obtained  by the same arguments used for \eqref{eqKINF} under {\bf (H2')}. 	
\item[(iii)]  Given $(K,\Lambda)$, we consider the equation for $(Y,Z^Y)$. This is a mean-field linear BSDE on infinite-horizon
 \beq 
dY_t &=&  \Big( \vartheta_t + G\trans(Y_t - \eee[Y_t]) + \hat G\trans \eee[Y_t] + J \trans(Z_t^Y - \eee[Z_t^Y]) + \;  \hat{J}\trans\eee[Z_t^Y] \Big) dt  \nonumber \\
& & \hspace{2cm} 	+ \; Z_t^Y dW_t,  \;\;\;\;\; t \geq 0	\label{eqYinfi}
 \enq
where $G$ $\hat G$, $J$, $\hat J$ are constant coefficients in $\R^{d\times d}$, and  $\vartheta$ is a random process in $L^2_{\F}(\Omega\times\R_+,\R^d)$  defined by
\begin{align*}
	& G \; = \;  \rho \,\, \I_d  - B + C S^{-1} U, 
	\\
	& \hat G \; = \;  \rho \,\, \I_d - B - \tilde B + (C +\tilde C)\hat S^{-1} V,
	\\
	& J \; = \;  - D + F S^{-1} U,
	\\
	& \hat{J} \; = \;  - (D + \tilde D) + ( F + \tilde  F) \hat S^{-1} V, 
	\\
	&\vartheta_t \; = \;  - M_t - K (\beta_t - \eee[\beta_t] ) - \Lambda \eee[\beta_t] - D \trans K (\gamma_t - \eee[\gamma_t]) - \hat D \trans K \eee[\gamma_t]
	\\
	& \qquad  \;\;\; + U \trans S^{-1} \big( H_t - \eee[H_t] + F \trans K (\gamma_t - \eee[\gamma_t])\big)  + V \trans \hat S^{-1} \big(\eee[H_t] + (F + \tilde F)\trans K \eee[\gamma_t] \big),
\end{align*}
with
\begin{equation} \label{Sinfi}
\left\{
\begin{array}{rcl}
S  &=&  N + F \trans K F \; = \; S(K), \\
\hat S &=& N + \tilde N + (F + \tilde F)\trans K (F + \tilde F)  \; = \; \hat S(K), \\
U &=& I + F\trans K D + C\trans K \; = \; U(K), \\
V  &=& I + \tilde I + (F +\tilde F)\trans K (D +\tilde D) +   (C + \tilde C)\trans \Lambda \; = \; V(K,\Lambda). 
\end{array}
\right.
\end{equation}
Although in many practical case an explicit solution is possible (see below), there are no general existence results for such a mean-field BSDE on infinite-horizon, to the best of our knowledge. We then assume what follows.

\begin{itemize}
\item[\textbf{(H4')}] There exists a solution $(Y,Z^Y)$ $\in$  $L^2_{\F}(\Omega\times \R_+,\R^d)\times L^2_{\F}(\Omega\times \R_+,\R^d)$ to \eqref{eqYinfi}. 
\end{itemize}

\vspace{1mm}

\begin{remark}
\label{remH4}
{\rm
In many practical cases, {\bf (H4')}  is satisfied and explicit expressions for $Y$ may be available. We list here  some examples.
\begin{itemize}
\item[-] In the case where $\beta=\gamma=H=M\equiv0$, so that $\vartheta$ $\equiv$ $0$, we see that $Y$  $=$ $Z^Y$ $\equiv 0$ is a solution to \eqref{eqYinfi}, and {\bf (H4')} clearly holds. 
\item[-] If $\beta,\gamma,H,M$ are deterministic (hence, all the coefficients are non-random), the process $Y$ is deterministic as well, that is $Z^Y \equiv 0$ and $\eee[Y]=Y$. Then, the mean-field BSDE becomes a standard linear ODE:
\beqs
dY_t &=& \big( \vartheta_t + \hat G\trans Y_t \big) dt, \;\;\; t \geq 0. 
\enqs
In the one-dimensional case $d=1$, we get
\begin{equation*}
Y_t = - \int_t^\infty e^{- \hat G (s-t)} \vartheta_s ds, \;\;\; t \geq 0.
\end{equation*}
Therefore, when $\hat G - \rho > 0$,  i.e., $(C +\tilde C)\hat S^{-1} V$ $>$ $B+\tilde B$, we have by the Jensen inequality and the Fubini theorem
\begin{equation*}
\int_0^\infty \!\! e^{-\rho t} Y_t^2 dt \leq \tilde c_1 \int_0^\infty \!\! \int_t^\infty \!\! e^{-\rho t}  e^{-\hat G (s-t)} \vartheta_s^2 ds \, dt \leq \tilde c_2 \int_0^\infty \!\! e^{-\rho s} \vartheta_s^2 ds < \infty,
\end{equation*}
for suitable constants $\tilde c_1, \tilde c_2>0$, so that {\bf (H4')} is satisfied. In the multi-dimensional case $d>1$, if $\beta, \gamma,H,M$ are constant vectors (hence, $\vartheta$ is constant as well), we have $Y_t = Y$, with
\begin{equation*}
Y = - (\hat G^{-1}) \trans \vartheta,
\end{equation*}
and {\bf (H4')} is clearly satisfied. 

\item[-] In many relevant cases, the sources of randomness of the state variable and the coefficients in the payoff are independent. More precisely, let us consider a problem with two independent Brownian motions $(W^1,W^2)$ (to adapt the formulas above, we proceed as in Remark \ref{remWmulti}). Assume that only $W^1$ appears in the controlled mean-field  SDE:
\beqs
dX^\alpha_t &=&  \big(\beta_t + B X^\alpha_t + \tilde B \eee[X^\alpha_t] + C \alpha_t + \tilde C \eee[\alpha_t]\big) dt 	\\
& & \;\;\;	+ \; \big(\gamma^1_t + D^1 X^\alpha_t + \tilde D^1 \eee[X^\alpha_t] + F^1 \alpha_t + \tilde F^1 \eee[\alpha_t]\big) dW_t^1,
\enqs
where $\beta,\gamma^1$ are deterministic processes. On the other hand,  the coefficients $M,H$ in the payoff are adapted to the filtration of $W^2$ and independent from $W^1$. In this case, the equation for $\big(Y,Z^Y=(Z^{1,Y},Z^{2,Y}) \big)$ writes
\beqs 
dY_t \!\! &=&  \!\! \Big( \vartheta_t + G\trans(Y_t - \eee[Y_t]) + \hat G\trans \eee[Y_t] + (J^1)\trans(Z_t^{1,Y} - \eee[Z_t^{1,Y}]) + \;  (\hat{J^1})\trans\eee[Z_t^{1,Y}] \Big) dt  \nonumber \\
& & \hspace{2cm} 	+ \; Z_t^{1,Y} dW^1_t + Z_t^{2,Y} dW^2_t,  \;\;\;\;\; t \geq 0,	
 \enqs
where we notice that $Z^{2,Y}$ does not appear in the drift as the corresponding coefficients are zero. Notice that  the process $(\vartheta_t)_t $ is adapted with respect to the filtration of $W^2$, while the other coefficients are constant.  It is then natural to look for a solution $Y$ also adapted w.r.t.~the filtration of $W^2$, hence for which $Z^{1,Y}$ $\equiv$ $0$, leading to the equation: 
\beqs
dY_t &=&  \Big( \vartheta_t + G\trans(Y_t - \eee[Y_t]) + \hat G\trans \eee[Y_t]  \Big) dt  + Z_t^{2,Y} dW^2_t, \;\;\; t \geq 0.
\enqs
For simplicity, let us consider the one-dimensional case $d=1$. Taking expectation in the above equation, we get a linear ODE for $\E[Y_t]$, and a linear BSDE for $Y_t - \eee[Y_t]$, given by
\beqs
d \eee[Y_t]&=&  \big( \E[\vartheta_t] + \hat G \eee[Y_t] \big) dt, \;\;\; t \geq  0,  \\
d(Y_t - \eee[Y_t]) &=& \Big( \vartheta_t - \E[\vartheta_t] + G (Y_t - \eee[Y_t]) \Big) dt +   Z_t^{2,Y} dW^2_t, \;\;\; t \geq 0,
\enqs
which lead to
\begin{equation*}
Y_t  \; = \;  - \int_t^\infty e^{-G(s-t)} \vartheta_s ds - \int_t^\infty \Big(e^{-\hat G(s-t)} - e^{- G(s-t)}\Big) \eee[\vartheta_s] ds.
\end{equation*}
Provided that $G -\rho, \hat G - \rho>0$, and recalling that $\vartheta$  $\in$ $L^2_{\F}(\Omega\times\R_+,\R^d)$, condition {\bf (H4')} is satisfied by the same estimates as above. See \cite{AidBaseiPham} and Section \ref{Sec:ExGoods} for practical examples from mathematical finance with such properties. \qed
\end{itemize}}
\end{remark}

\item[(iv)]  
Given $(K,\Lambda,Y,Z^Y)$ the equation for $R$ is a linear ODE, whose unique solution is explicitly given  by 
\beq \label{eqRinfi}
R_t &=& \int_t^\infty e^{-\rho (s-t)} h_s ds, \;\;\; t \geq 0, 
\enq
where the deterministic function $h$ is defined, for $t$ $\in$ $\R_+$, by
\beqs
h_t &=&  \eee\big[ - \gamma_t\trans K \gamma_t - \beta_t\trans Y_t - 2\gamma_t\trans Z_t^Y  + \xi_t\trans S^{-1} \xi_t \big] 
- \eee[\xi_t]\trans S^{-1} \eee[\xi_t] + O_t\trans \hat S^{-1} O_t,
\enqs
with 
\begin{equation} \label{xiinfi}
\left\{
\begin{array}{rcl}
\xi_t &=& H_t + F \trans K \gamma_t  + C \trans Y_t + F\trans Z_t^Y, \\
O_t &=& \E[H_t] + (F + \tilde F)\trans K \E[ \gamma_t] + (C + \tilde C)\trans\E[Y_t] + ( F+ \tilde F)\trans \E[Z_t^Y].
\end{array}
\right.
\end{equation}
Under assumptions {\bf (H1')}(i), {\bf (H2')}(ii) and {\bf (H4')},  we see that $\int_0^\infty e^{-\rho t} |h_t| dt$ $<$ $\infty$, from which we obtain that 
$R_t$ is well-defined for all $t$ $\geq$ $0$. Therefore,  $e^{-\rho t} |R_t|$ $\leq$ $\int_t^\infty e^{-\rho s} |h_s| ds$ $\rightarrow$ $0$ as $t$ goes to infinity. 
\end{itemize}

\paragraph{Final step.}  
Given $(K,\Lambda,Y,Z^Y,R)$ solution to \reff{eqKINF}, \reff{eqLINF}, \reff{eqYinfi}, \reff{eqRinfi},  hence to the system \reff{sysKinfi},  the function 
\beqs
w_t(x,\bar x) &=& (x - \bar x)\trans K (x - \bar x) + \bar x\trans \Lambda \bar x + 2 Y_t\trans x + R_t,
\enqs
satisfies the growth condition \reff{growthwinfi}, and by construction  the conditions (i)-(ii) of the verification theorem in Lemma \ref{prop:verifINF}. 
Let us now  consider as in the finite-horizon case the candidate for the optimal control given by 
\beq
\alpha_t^* &=& a_t^0(X_t^*,\E[X_t^*]) + a_t^1(\E[X_t^*]) \nonumber\\
&=& - S^{-1} U (X_t-\E[X_t^*]) -  S^{-1} \big(\xi_t - \E[\xi_t] \big)  -   \hat S^{-1} \big(V  \E[X_t^*] +  O_t  \big), \;\; t \geq 0,  \label{alphaoptinfi}
\enq
where $X^*$ $=$ $X^{\alpha^*}$ is the state process with the feedback control $a_t^0(X_t^*,\E[X_t^*]) + a_t^1(\E[X_t^*])$, and $S,\hat S,U,V$, $\xi,O$ are recalled in \reff{Sinfi}, \reff{xiinfi}. 
With respect to the finite-horizon case in Proposition \ref{profini}, the main technical issue is to check that $\alpha^*$ satisfies  the admissibility condition in $\Ac$. We need to make an additional condition on the discount factor:

\vspace{5mm}
\noindent \textbf{(H5')} \hspace{0.8cm} $\rho >  2 \max \Big\{ |B-CS^{-1}U| + |D-FS^{-1}U|^2, \,\, |(B+\tilde B)-(C+\tilde C) \hat S^{-1} V| \Big\}$.
\vspace{5mm}

\noindent From the  expression of $\alpha^*$ in \eqref{alphaoptinfi}, we see that $X^*=X^{\alpha^*}$ satisfies, 
\beqs
dX^*_t &=& b^*_t dt + \sigma^*_t dW_t, \;\;\; t \geq 0, 
\enqs
with
\beqs
b^*_t \; = \;  \beta^*_t + B^* (X^*_t - \eee[X^*_t]) + \tilde B^* \eee[X^*_t], & & 
\sigma^*_t  \;= \;   \gamma^*_t + D^* (X^*_t - \eee[X^*_t]) + \tilde D^* \eee[X^*_t],
\enqs
where we set
\begin{gather*}
B^* \; = \;  B-CS^{-1} U,
\qquad\qquad
\tilde B^* = (B + \tilde B) - (C+\tilde C) \hat S^{-1} V,
\\
D^* \; = \;  D - FS^{-1} U ,
\qquad\qquad
\tilde D^* \; = \;  (D + \tilde D) - (F + \tilde F) \hat S^{-1}V,
\\
\beta^*_t \; = \;  \beta_t  - CS^{-1} (\xi_t-\eee[\xi_t]) - (C + \tilde C) \hat S^{-1} O_t,
\\
\gamma^*_t \;= \;  \gamma_t  - FS^{-1} (\xi_t-\eee[\xi_t]) - (F + \tilde F )\hat S^{-1} O_t.
\end{gather*}
In order to prove the admissibility condition in $\Ac$:  $\int_0^\infty e^{-\rho t} \eee[|\alpha^*_t|^2] dt$ $<$ $\infty$, and as $\xi$ $\in$ $L_\F^2(\Omega\times\R_+,\R^m)$, 
$O$ $\in$ $L^2(\R_+,\R^m)$, it suffices to verify that 
\begin{equation}
\label{eqadm2}
\int_0^\infty e^{-\rho t} \eee\big[\big|X^*_t - \eee[X^*_t]\big|^2\big] dt < \infty,
\qquad\quad
\int_0^\infty e^{-\rho t} \big| \eee[X^*_t] \big|^2 dt < \infty.
\end{equation}
Let us start with the second integral in \eqref{eqadm2}. By the It\^o's formula and the Young inequality, for each $\eps>0$ we have (using shortened bar notations, see Remark \ref{rem:notations}, e.g. $\bar X^*$ $=$ $\E[X^*]$)
\beq
\frac{d}{dt}e^{-\rho t}|\bar X_t^*|^2 &=&  e^{-\rho t}\big(-\rho |\bar X_t^*|^2 + 2 (\bar b_t^*)\trans  \bar X_t^* \big) \nonumber \\
	&\leq& e^{-\rho t}\Big(-\rho |\bar X_t^*|^2 + 2 |\bar\beta_t^*| |\bar X_t^*| + 2(\bar X_t^*)\trans\tilde B^*\bar X_t^* \Big) \label{estimB*} \\
	&\leq & e^{-\rho t}\Big[ \Big( -\rho + 2|\tilde B^*| + \eps \Big) |\bar X_t^*|^2 + c_\eps \big( |\bar\beta_t|^2 +|O_t|^2 \big) \Big], \nonumber 
\enq
where $c_\eps>0$ is a suitable constant. We now set 
\begin{equation*}
\zeta^* \; = \;  \big|\eee[X_0]\big|^2 + c_\eps  \int_0^\infty e^{-\rho t} \big[|\bar\beta_t|^2 + |O_t|^2 \big] dt, 
\qquad
\eta^* \; = \;  -\rho + 2|B^*+\tilde B^*| + \eps,
\end{equation*}
and notice that $\zeta^*$ $<$ $\infty$ by {\bf (H1')}-{\bf (H2')}, while $\eta^* <0$ for $\eps$ small enough, by {\bf (H5')}. 
By applying the Gronwall inequality, we then get
\begin{equation*}
\int_0^\infty e^{-\rho t}\big|\eee[X^*_t]\big|^2 dt  \; \leq \;  \zeta^* \int_0^\infty e^{(\eta^*) t} dt \; < \; \infty. 
\end{equation*}	
For the first integral in \eqref{eqadm2}, by similar estimates we have
\beq
& &\frac{d}{dt}\eee\big[e^{-\rho t}|X_t^* - \bar X_t^*|^2\big] \nonumber  \\
&=& e^{-\rho t}\eee\big[-\rho |X_t^* - \bar X_t^*|^2 + 2( b_t^* - \bar b_t^*)\trans(X_t^* -  \bar X_t^*) + |\sigma_t^{*}|^2 \big] \nonumber \\
&\leq& e^{-\rho t}\eee\Big[-\rho |X_t^* - \bar X_t^*|^2 + 2 \Big( |\beta_t^* - \bar \beta_t^*| |X_t^* -  \bar X_t^*| + (X_t^* -  \bar X_t^*)\trans B^*( X_t^* -  \bar X_t^*) \Big) \nonumber  \\
& &\qquad\qquad + \;  2 \Big( |\gamma_t^*|^2 + |D^*|^2|X_t^* - \bar X_t^*|^2 + |\tilde D^*|^2 |\bar X_t^*|^2 \Big) \Big] \label{estimtildeB*}  \\
&\leq& e^{-\rho t}\eee\Big[ \Big( -\rho + 2|B^*| + 2|D^*|^2 + \eps \Big) |X_t^* - \bar X_t^*|^2  +  
 c_\eps\big( |\beta_t|^2 + |\gamma_t|^2 + |O_t|^2 + |\xi_t|^2 + |\bar X_t^*|^2 \big) \Big], \nonumber
\enq
for a suitable $c_\eps>0$. Recall that $\int_0^\infty e^{-\rho t} |\eee[X^*_t]|^2 dt < \infty$, we deduce as  above  by Gronwall inequality, and under {\bf (H5')}, that 
\begin{equation*}
\int_0^\infty e^{-\rho t} \eee\big[\big|X^*_t - \eee[X^*_t]\big|^2\big] dt \; <\; \infty,
\end{equation*}	
which shows that $\alpha^*$ $\in$ $\Ac$. 


\vspace{1mm}

\begin{remark}
\label{RemSimplH5}
{\rm In some specific cases, {\bf (H5')} can be weakened. For example, assume that 
\beqs
B-CS^{-1} U \;\leq \; 0, \;\;\;  (B + \tilde B) - (C+\tilde C) \hat S^{-1} V \; \leq \;  0, \;\;  \gamma_t = \tilde D = F = \tilde F =0.
\enqs
In this case, the matrices $B^* ,\tilde B^* $ are negative definite, and the estimates in \reff{estimB*}, \reff{estimtildeB*} are simpler and some terms disappear; also, the coefficient $2$ is no longer necessary (see Remark \ref{RemSimplH3}). Assumption {\bf (H5')} can be simplified into
\begin{equation*}
\pushQED{\qed} 
\rho \;>\; |D|^2. \qedhere 
\popQED
\end{equation*}
}
\end{remark}


To sum up the arguments of this section, we have proved the following result.

\begin{Theorem} \label{thm:optimalINF}
	Under assumptions  {\bf (H1')}-{\bf (H5')}, the optimal control for the LQMKV  pro\-blem on infinite-horizon  \eqref{pb:payoffINFI} is given  by  \reff{alphaoptinfi}, and  
	the corresponding value of the problem is 
	\beqs
	V_0 &=& J(\alpha^*) \; = \;  \eee\big[ (X_0- \eee[X_0])\trans K (X_0- \eee[X_0]) \big] +  \eee[X_0]\trans \Lambda  \eee[X_0] + 2 \eee\big[Y_0\trans X_0] + R_0.
	\enqs
\end{Theorem}

\begin{remark}  \label{remmultiinfi}
{\rm
The remarks in Section \ref{Sec:remarks} can be immediately adapted to the infinite-horizon framework. In particular, as in Remark \ref{remH22}, 
one can have existence of a solution to \reff{eqKINF}-\reff{eqLINF}, even when condition {\bf (H2')}  is not satisfied, and obtain the optimal control as in \reff{alphaoptinfi} provided that {\bf (H4')}-{\bf (H5')} are satisfied. 
On the other hand,  the model considered here easily extends to the case where several independent Brownian motions are present, as described in 
Remark \ref{remWmulti}.
}
\qed
\end{remark}

\section{Applications}
\label{Sec:applications}

\subsection{Portfolio liquidation with trade crowding}

We focus on the now classical problem of portfolio liquidation, initiated in \cite{almcri00}, but taking also into account the permanent price impact due to similar participants as in \cite{carleh16}.  We consider a finite-horizon model and compute the optimal control by Theorem \ref{thm:optimal}.

A trader has to liquidate a certain number of shares $x_0 \geq 0$ before a terminal time $T>0$.  She controls the trading speed $\alpha_t$ through time, and her inventory $X_t^\alpha$ evolves according to
\beq \label{Xinv}
X^\alpha_t &=& x_0 + \int_0^t  \alpha_s ds,
\enq
for $0 \leq t \leq T$ (notice that for a liquidation problem, i.e., $x_0$ $>$ $0$, the control $\alpha$ will be mostly negative corresponding to a sale of the stock). The price $S=S^\alpha$ of the asset is subject to a permanent market impact generated by the trading of other identical market participants, who act cooperatively with each other. This differs from the mean-field game situation considered in \cite{carleh16}, where the investors act strategically in a Nash equilibrium. Here, we focus on a Pareto optimality framework, and in the asymptotic regime of large number of investors, the stock price is given by
\beqs
S_t^\alpha &=& S_t^0 + \nu \int_0^t \E[\alpha_s] ds,
\enqs
where $S^0$ is an exogenous $\F$-adapted process representing the asset price in absence of trading, and $\nu$ $\geq$ $0$ is a positive parameter modeling the linear permanent impact due to the average trading of all market participants. The objective of the investor is then to minimize, over her trading speed $\alpha$, the expected total liquidation cost
\beqs
J(\alpha) &=& \E\Big[ \int_0^T \Big( q (X^\alpha_t)^2 + \alpha_t( S^\alpha_t + \eta \alpha_t) \Big) dt + p (X^\alpha_T)^2 \Big],
\enqs
where $\eta$ $>$ $0$ is a constant parameter modelling the linear temporary market impact, $q$ $\geq$ $0$ is a penalty parameter on the current inventory, and $p$ $\geq$ $0$ is a parameter penalizing the remaining inventory at maturity $T$. Since $S^\alpha_t = S^0_t + \nu (\eee[X^\alpha_t] - x_0)$, the cost functional can be rewritten as
\beqs
J(\alpha) &=& \E\Big[ \int_0^T \Big( q (X^\alpha_t)^2 + \eta \alpha^2_t + ( S_t^0 - \nu x_0) \alpha_t + \nu\eee[\alpha_t]\eee[X^\alpha_t] \Big) dt + p (X^\alpha_T)^2 \Big].
\enqs
W.l.o.g., we can assume that $\eta$ $=$ $1$ by normalizing the other coefficients: $q$ $\rightarrow$ $q/\eta$, $S_t^0$ $\rightarrow$ $S_t^0/\eta$, 
$\nu$ $\rightarrow$ $\nu/\eta$, and  $p$ $\rightarrow$ $p/\eta$. This problem fits into the framework in Sections \ref{Sec:intro} and \ref{Sec:assumpt}, with $d=m=1$ (one-dimensional state variable and control), and coefficients, in the notations of \eqref{pb:coeffSDE} and \eqref{pb:coeffPayoff},   given for $t \in [0,T]$ by
\begin{equation*}
C_t = 1, \qquad
Q_t = q, \qquad
N_t = 1, \qquad
H_t = \frac{S_t^0 - \nu x_0}{2}, \qquad
\tilde I_t = \nu,    \qquad
P_t =  p.
\end{equation*}
The other coefficients, including the discount rate, are identically zero. Notice that the coefficient $H_t$ is stochastic.  Conditions {\bf (H1)}-{\bf (H2)} are satisfied once
\beqs
\int_0^T \eee[|S^0_t|^2] dt \; < \;  \infty, & \mbox{ and } & q \; \geq \; \nu^2. 
\enqs
The  case where $q$ $=$ $0$, even when $\nu$ $>$ $0$ (hence when {\bf (H2)} is not satisfied), 
can also be dealt with, as pointed out in Remark \ref{remH22}, and discussed more precisely in Remark \ref{remqzero}.  
From Theorem \ref{thm:optimal}, the optimal control is given by
\beq \label{alphainter}
\alpha_t^* &=& - K_t(X_t^* - \E[X_t^*]) - (\Lambda_t + \nu) \E[X_t^*] - Y_t  - \frac{S_t^\nu}{2}, \;\;\; 0 \leq t\leq T,
\enq
where we set $S_t^\nu$ $:=$ $S_t^0 - \nu x_0$.  Here, the Riccati equation for $K$ writes as 
\begin{equation} \label{RicK}
\begin{cases}
\dot K_t + q - K^2_t \; = \;  0, \quad t \in [0,T],  \\
K_T \; = \;  p,
\end{cases}
\end{equation}
whose solution is explicitly given by 
 \beqs
K_t &=& \sqrt{q}  \frac{ \sqrt{q}\sinh(\sqrt{q}(T-t)) + p\cosh(\sqrt{q}(T-t))}
{ p  \sinh(\sqrt{q}(T-t)) + \sqrt{q}\cosh(\sqrt{q}(T-t))}, \;\;\; 
0 \leq t \leq T. 
\enqs 
The Riccati equation for $\Lambda$ writes as  
\begin{equation} \label{RicL}
\begin{cases}
\dot \Lambda_t + q - (\Lambda_t + \nu)^2 \; = \;  0,   \quad t \in [0,T],\\
\Lambda_T \; = \;  p.
\end{cases}
\end{equation}
By writing $K^\nu$ $:=$ $\Lambda + \nu$, we get the same equation as above for $K$, with final condition $p+\nu$, and so
\beq \label{Knu}
K_t^\nu  &=& \sqrt{q}  \frac{ \sqrt{q}\sinh(\sqrt{q}(T-t)) + (p+\nu)\cosh(\sqrt{q}(T-t))}
{(p+\nu) \sinh(\sqrt{q}(T-t)) + \sqrt{q}\cosh(\sqrt{q}(T-t))}, \;\;\; 
0 \leq t \leq T,
\enq
while we notice that $K$ $=$ $K^0$ for $\nu$ $=$ $0$. 
The equation for $Y$ reads as 
\begin{equation} \label{dynYinter}
\begin{cases}
dY_t \; = \;  \big[ \frac{K_t(S^\nu_t - \eee[S^\nu_t])+  K_t^\nu\eee[S^\nu_t]}{2} +  K_t (Y_t - \eee[Y_t]) + K_t^\nu \eee[Y_t] \big] dt 
+ Z^Y_t dW_t, \quad t \geq 0,
\\
Y_T \; = \; 0.
\end{cases}
\end{equation}
By taking expectation in the above equation, and then writing  $dY_t = d(Y_t - \eee[Y_t]) + d\eee[Y_t]$, we easily solve this linear mean-field BSDE:
\begin{equation} \label{Yinter}
\left\{
\begin{array}{rcl}
Y_t &=& - \int_t^T K_s e^{- \int_t^s K_udu} \frac{ (\eee[S^\nu_s|\mathcal{F}_t] - \eee[S^\nu_s]) }{2} ds + \E[Y_t], \\
\E[Y_t] &=& - \int_t^T K_s^\nu e^{- \int_t^s K_u^\nu  du} \frac{\eee[S^\nu_s]}{2} ds, \;\;\; 0 \leq t\leq T. 
\end{array}
\right.
\end{equation}
Moreover, from \reff{Xinv} and \reff{alphainter}, the mean of $X^*$ is governed by 
\beqs
d\E[X_t^*] &=& \E[\alpha_t^*] dt \; = \;  - \Big( K_t^\nu \E[X_t^*] +  \E[Y_t]   + \frac{\E[S_t^\nu]}{2} \Big) dt, \;\;\; \E[X_0^*] \; = \; x_0,
\enqs
hence explicitly given by 
\beq 
\E[X_t^*] &=& x_0 e^{-\int_0^t K_u^\nu du}  - \int_0^t e^{-\int_s^t K_u^\nu  du} \Big( \E[Y_s]  + \frac{\E[S_s^\nu]}{2}  \Big) ds \nonumber \\
&=&   x_0 e^{-\int_0^t K_u^\nu du} + \frac{1}{2} \int_0^t e^{-\int_s^t K_u^\nu  du}  
\Big(   \int_s^T K_r^\nu e^{-\int_s^r K_u^\nu du} \E[S_r^\nu] dr  - \E[S_s^\nu] \Big) ds,  \label{meanXinter}
\enq
by \reff{Yinter}. 
By integrating the function $K^\nu$, we have 
\beqs
e^{-\int_t^s K_u^\nu du} 
&=& \frac{K_t^\nu}{K_s^\nu} \frac{\omega_\nu(T-t)}{\omega_\nu(T-s)} \; = \; \frac{\pi_\nu(T-t)}{\pi_\nu(T-s)},
\enqs
where we set for $\tau$ $\in$ $[0,T]$, 
\begin{equation} \label{omegapi}
\left\{
\begin{array}{ccc}
\omega_\nu(\tau) &=& \frac{p+\nu}{\sqrt{q} \sinh(\sqrt{q}\tau) + (p+\nu)\cosh(\sqrt{q}\tau)} \;\;  \in \; (0,1], \\
\pi_\nu(\tau) &=& \frac{\sqrt{q}}{\sqrt{q} \cosh(\sqrt{q}\tau) + (p+\nu)\sinh(\sqrt{q}\tau)} \;\;  \in \; (0,1].   
\end{array}
\right.
\end{equation}
By plugging into \reff{Yinter}, \reff{meanXinter}, and recalling that $S^\nu$ $=$ $S^0-\nu x_0$, 
the optimal control in \reff{alphainter} is then expressed as
\beq
\alpha_t^* 
&=&    - K_t X_t^*  +  \frac{1}{2} \Big(  \int_t^T K_t  \frac{\omega_0(T-t)}{\omega_0(T-s)} \E[S^0_s|\mathcal{F}_t] ds - S_t^0 \Big)  \;  + \;  x_0 \frac{\nu}{2} \pi_\nu(T-t) 
\nonumber \\
& &  \;\;\; - \;  x_0 (K_t^\nu - K_t) \Big\{ \frac{\pi_\nu(T)}{\pi_\nu(T-t)}
+  \frac{\nu}{2} \int_0^t  \frac{\pi_\nu(T-s)}{\pi_\nu(T-t)} \pi_\nu(T-s) ds \Big\} \nonumber \\
& & \;\;\; - \; \frac{K_t^\nu - K_t}{2} \int_0^t  \frac{\pi_\nu(T-s)}{\pi_\nu(T-t)}
\Big(   \int_s^T  K_s^\nu \frac{\omega_\nu(T-s)}{\omega_\nu(T-r)} \E[S_r^0] dr  - \E[S_s^0] \Big) ds \nonumber \\
& & \;\;\; + \; \frac{1}{2} \int_t^T \Big( K_t^\nu  \frac{\omega_\nu(T-t)}{\omega_\nu(T-s)}  - K_t  \frac{\omega_0(T-t)}{\omega_0(T-s)}  \Big) \E[S_s^0] ds,  \;\;\;\;\;  0 \leq t \leq T.  \label{alphaoptexpli1}
\enq

\vspace{1mm}

\begin{figure}[h!]
	\begin{minipage}{0.49\textwidth}
		\centering
		\includegraphics[width=\textwidth]{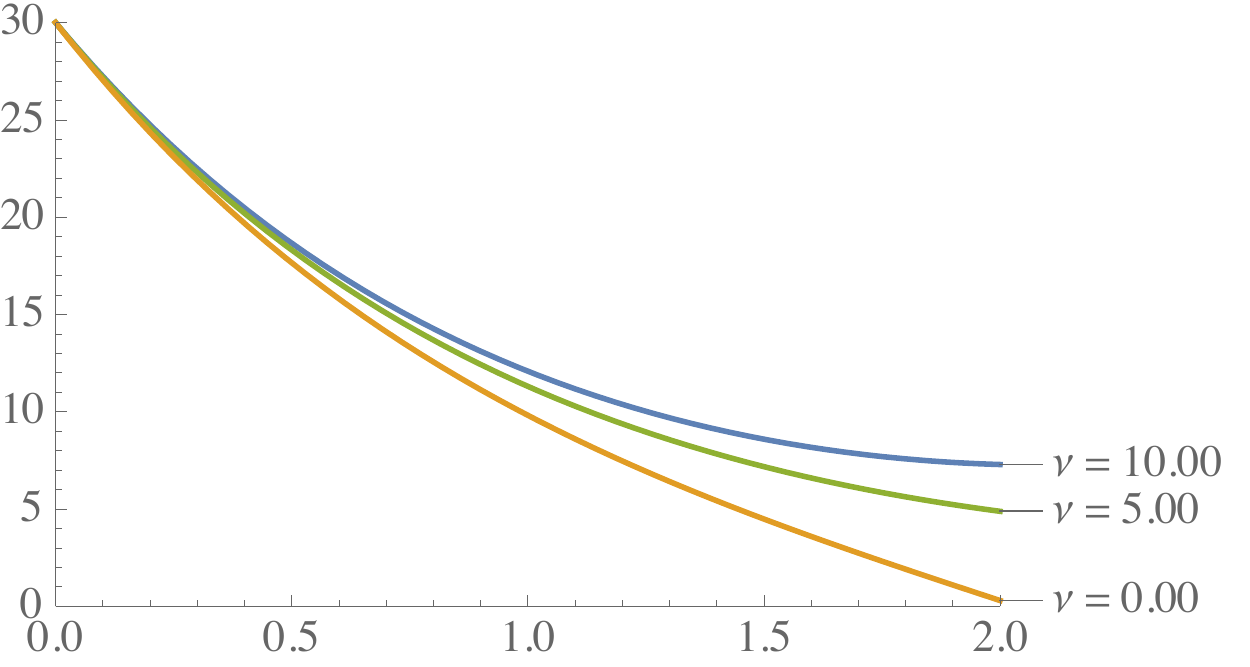}
		\caption{\footnotesize{Graph of  $t$ $\mapsto$ $\E[X_t^*]$, for $p$ $=$ $10$, $q$ $=$ $1$, $x_0$ $=$ $30$, $S_0$ $=$ $10$,  $T=2$, and by varying $\nu$.}}
		\label{fig:evolutionmeaninventory1}
	\end{minipage}
	\hfill 
	\begin{minipage}{0.49\textwidth}
		\centering
		\includegraphics[width=\textwidth]{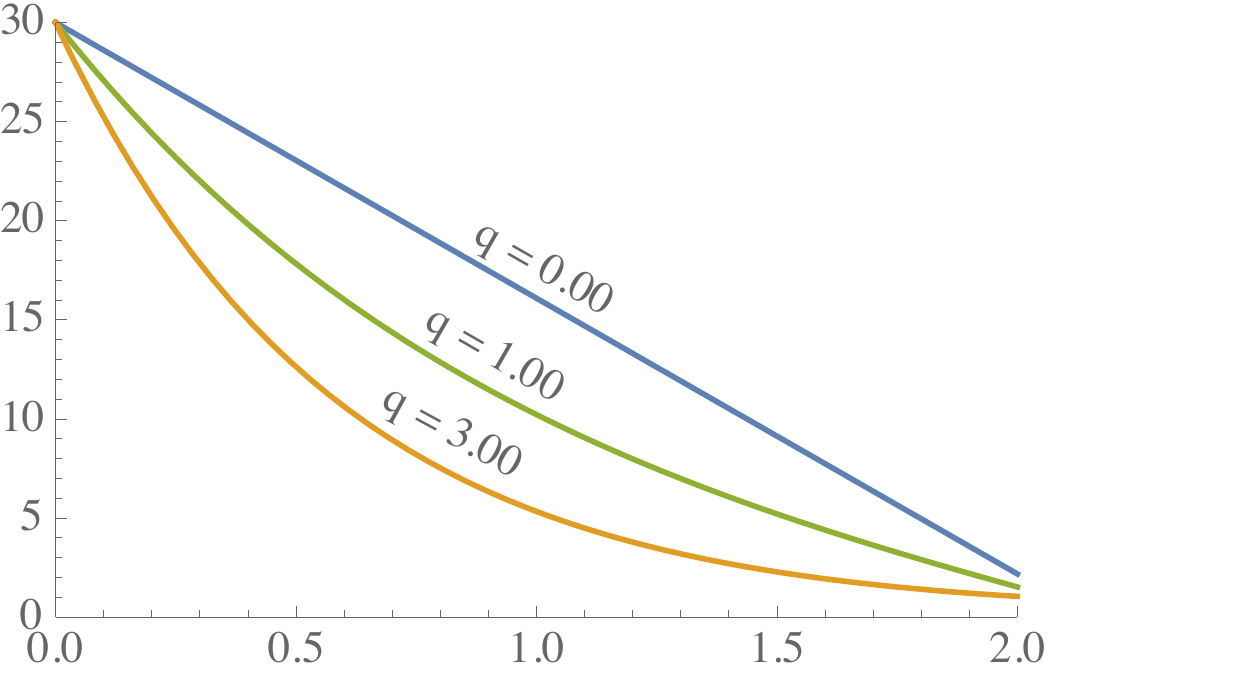}
		\caption{\footnotesize{Graph of  $t$ $\mapsto$ $\E[X_t^*]$, for $p$ $=$ $10$, $\nu$ $=$ $1$, $x_0$ $=$ $30$, $S_0$ $=$ $10$,  $T=2$, and by varying $q$.}}
		\label{fig:evolutionmeaninventory2}
	\end{minipage}
\end{figure}

\begin{remark}
{\rm  In absence of permanent price impact, i.e., $\nu$ $=$ $0$, the optimal execution strategy is given by
\beqs
\alpha_t^{0,*} &=& - K_t X_t^* + \frac{1}{2} \Big( \int_t^T K_t  \frac{\omega_0(T-t)}{\omega_0(T-s)} \E[S^0_s|\mathcal{F}_t] ds - S_t^0 \Big), \\
& =:& - K_t X_t^*  + \alpha_t^{0,S},  \;\;\;\;\;  0 \leq t \leq T.
\enqs
The first term prescribe the agent to liquidate towards the zero inventory target, while the second term $\alpha^{0,S}$ is an incentive to buy or sell depending on whether the weighted average of expected future value of the stock is larger or smaller than its current value. In particular, when $S^0$ is a martingale, we have
\beqs
\alpha_t^{0,S} &=& - \frac{S_t^0}{2} \pi_0(T-t),
\enqs
which is nonpositive (whenever $S^0$ is nonnegative), meaning a sale of the asset.  The four additional terms in \reff{alphaoptexpli1} are deterministic, and involve the initial capital $x_0$,  the expected value of the asset price $S^0$, and measure quantitatively the impact of the permanent price impact in the optimal execution strategy where we notice (e.g.~by maximum principle for the ODE \reff{RicL}) that $K_t^\nu$ $\geq$ $K_t$.  
In the case where $S^0$ is a martingale, we have
\beqs
\alpha_t^* - \alpha_t^{0,*} &=&      x_0 \Big\{ \frac{\nu}{2} \pi_\nu(T-t)  -    (K_t^\nu - K_t) \Big[ \frac{\pi_\nu(T)}{\pi_\nu(T-t)}
+  \frac{\nu}{2} \int_0^t   \frac{\pi_\nu(T-s)}{\pi_\nu(T-t)} \pi_\nu(T-s) ds \Big] \Big\} \\
& &  \; + \; \frac{S_0^0}{2} \Big\{ (K_t^\nu - K_t) \int_0^t   \frac{\pi_\nu(T-s)}{\pi_\nu(T-t)} \pi_\nu(T-s) ds  +   \big( \pi_0(T-t) - \pi_\nu(T-t) \big)
\Big\} .
\enqs
On the other hand, by applying It\^o's formula to \reff{alphainter}, and using \reff{RicK}-\reff{RicL}-\reff{dynYinter}, while recalling  the dynamics \reff{Xinv} of the inventory, we have 
\beqs
d\big(\alpha_t^* + \frac{S_t^0}{2}\big)  \; = \; d\big(\alpha_t^* + \frac{S_t^\nu}{2}\big)  &=& q X_t^* dt   - Z_t^Y dW_t,
\enqs
which shows the notable robust property on the optimal execution strategy
\beqs
\alpha_t^* + \frac{S_t^0}{2} - q \int_0^t X_s^* ds, \;\; 0 \leq t \leq T,  &  & \mbox{ is a martingale.}
\enqs
In particular, when $S^0$ is a martingale, this implies that: 
\beqs
\E[\alpha_t^*] &=& \alpha_0^* + q \int_0^t \E[X_s^*] ds, \;\;\; 0 \leq t \leq T.  
\enqs
As $\E[X_t^*]$ $=$ $x_0 + \int_0^t \E[\alpha_s^*] ds$, we see that $E(t)$ $:=$ $\E[X_t^*]$ satisfies the second-order ODE: $E''(t)$ $=$ $q E(t)$, with initial condition: 
$E(0)$ $=$ $x_0$, $E'(0)$ $=$ $\alpha_0^*$,  from which we deduce that 
\beq
E(t) \; := \; \E[X_t^*] &=& x_0 \cosh(\sqrt{q}t) + \frac{\alpha_0^*}{\sqrt{q}} \sinh(\sqrt{q}t),  \;\;\; 0 \leq t \leq T, \label{Xoptmean} \\
&=& x_0 \Big[ \cosh(\sqrt{q}t) - \frac{K_0^\nu - \frac{\nu}{2}\pi_\nu(T)}{\sqrt{q}} \sinh(\sqrt{q}t) \Big] \; - \;  \frac{S_0^0 \pi_\nu(T)}{2\sqrt{q}} \sinh(\sqrt{q}t), \nonumber 
\enq
where we used expression \reff{alphaoptexpli1} of $\alpha_t^*$ for $t$ $=$ $0$.  Recalling the expressions of $K_0^\nu$ and $\pi_\nu(T)$ in \reff{Knu}, \reff{omegapi}, this gives an explicit formula for the optimal inventory on mean in terms of the parameters: 
$q$ (the risk aversion on the current inventory), $p$ (the penalty on the terminal inventory), $\nu$ (the mean-field permanent market impact), $x_0$ (the initial net position), $S_0^0$ (the initial value of the stock), and $T$ (the trading duration). 
It is easily seen from \reff{Xoptmean} that  $t$ $\mapsto$ $E(t)$ is convex.  
We plot in Figure \ref{fig:evolutionmeaninventory1} the evolution of the optimal inventory on mean 
$t$ $\rightarrow$ $E(t)$ for typical values of parameters: $p$ $=$ $10$, $q$ $=$ $1$, $x_0$ $=$ $30$, $S_0$ $=$ $10$, $T$ $=$ $2$, and 
by varying $\nu$.  We observe  that a smaller  mean-field parameter impact $\nu$ drives faster the  inventory to liquidation.  For $\nu$ $=$ $1$, we obtain $E(T)$ $=$ $1.52$. Finally, notice that as $p$ goes to infinity (corresponding to the null constraint on the terminal inventory $X_T$), then $E(t)$ converges to
\beqs
E_\infty(t) &=& x_0 \Big[  \cosh(\sqrt{q}t) - \frac{\cosh(\sqrt{q}T)}{\sinh(\sqrt{q}T)} \sinh(\sqrt{q}t) \Big], \;\;\; 0 \leq t \leq T,
\enqs
which is equal to zero for $t$ $=$ $T$, as expected, but  does not depend on the permanent impact parameter $\nu$, neither on $S_0^0$.  
}
\ep 
\end{remark}



\begin{remark} \label{remqzero}
{\rm  When $q$ $=$ $0$ (as in the original formulation of \cite{almcri00} with $\nu$ $=$ $0$), condition {\bf (H2)} is not satisfied when $\nu$ $>$ $0$. However, solutions  
$K$ and $\Lambda$ to the Riccati equations \reff{RicK} and \reff{RicL} still exist, given by $K_t$ $=$ $K_t^0$ and $\Lambda_t$ $=$ $K_t^\nu - \nu$, where $K_t^\nu$ is now defined by 
\beqs
K_t^\nu &=& \frac{p+\nu}{(p+\nu)(T-t) + 1}, \;\;\; 0 \leq t \leq T, 
\enqs
but $\Lambda$ may not be nonnegative in this case, whereas the condition $q$ $\geq$ $\nu^2$ in {\bf (H2)} ensures the nonnegativity of $\Lambda$.  Anyway, 
one can apply Proposition \ref{profini}, and the expression \reff{alphaoptexpli1} for the optimal control holds true. Some simplifications arise for $q$ $=$ $0$. Indeed, in this case 
$\omega_\nu$ $\equiv$ $1$, $\pi_\nu(T-t)$ $=$ $1/(1+(p+\nu)(T-t))$ $=$ $K_t^\nu/(p+\nu)$.  
Moreover, when $S^0$ is a martingale, the optimal mean inventory decreases linearly w.r.t. time, and is given by 
\beqs
\E[X_t^*] &=& x_0 \big[ 1 - \frac{\big(p + \frac{\nu}{2}\big)t}{1 + (p+\nu)T}  \big] - \frac{S_0^0}{2} \frac{t}{1 + (p+\nu)T}, \;\;\; 0 \leq t\leq T.
\enqs
We plot in Figure \ref{fig:evolutionmeaninventory2} the evolution of the optimal inventory on mean $t$ $\rightarrow$ $\E[X_t^*]$ for different values of  $q$, and observe that a larger  risk aversion parameter $q$ drives faster the inventory to liquidation. 
}
\ep
\end{remark}

\subsection{Production of an exhaustible resource}
\label{Sec:ExGoods}

We study an infinite-horizon model of substitutable production goods of exhaustible resource with a large number of producers,  inspired by  the papers 
\cite{guelaslio10} and \cite{chasir14}, see also \cite{gra16}.    
Let us consider a representative producer, and denote by $\alpha_t$ the quantity supplied at time $ t \geq 0$, and by $X^\alpha_t$ her current level of reserve in the good. 
As in \cite{guelaslio10}, we assume that the dyna\-mics of the reserve is stochastic with a noise proportional to the current level of reserves, hence evolving according to 
\beqs
dX^\alpha_t &=& - \alpha_t dt + \sigma X_t^\alpha dW^1_t,\;\;\;  t \geq 0,  \;\;\; X_0^\alpha \; = \; x_0, 
\enqs
where $\sigma$ $>$ $0$, and $W^1$ is a standard Brownian motion. The goal of each producer is the maximization of her profit, and we shall focus on a cooperative market equilibrium. Denote by $P^\alpha$  the price received by the producer when selling with policy $\alpha$. Following \cite{chasir14}, we assume a linear inverse demand rule, which models the price as a linear (decreasing) function of the quantity produced:
\beqs
P^\alpha_t &=& P_t^0 - \delta \alpha_t - \eps \E[\alpha_t], 
\enqs
where $\delta, \eps$ $>$ $0$ are a positive constants, and $P^0$ is some exogenous random process, adapted with respect to the filtration of some Brownian motion $W^2$ independent of $W^1$. The para\-meter $\eps$ measures the degree of interaction or the product substitutability, in the sense that the price received by an individual firm decreases as the other firms increase the supply of their goods. In particular, the case $\eps$ $=$ $0$ corresponds to independent goods. The objective of the representative producer is then formulated as the maximization over $\alpha$ of the gain functional\footnote{We thank Ren\'e Aid for insightful discussions on this example.}
\beqs
\E \Big[ \int_0^\infty e^{-\rho t} \big(  \alpha_t P^\alpha_t   -   \eta {\rm Var}(\alpha_t) - c \alpha_t( x_0 - X_t^\alpha) \big) dt \Big], 
\enqs
where $\rho > 0$ is the discount rate, 
$\eta$ $\geq$ $0$ is a nonnegative parameter penalizing high variation of the produced quantity, and the last term $\Cc(\alpha)$ $=$ 
$c \alpha_t( x_0 - X_t^\alpha)$, with $c$ $>$ $0$, represents the cost of extraction. In the beginning, this cost is negligible, and increases as the reserve is depleted. The control problem is  written equivalently as the minimization of the cost functional
\beqs
J(\alpha) &=& \E \Big[ \int_0^\infty e^{-\rho t} \Big( (\delta+\eta)(\alpha_t - \E[\alpha_t])^2 + (\delta+ \eps) \big(\E[\alpha_t]\big)^2  - c \alpha_t X_t^\alpha  
+ \alpha_t (cx_0 -  P_t^0)  \Big) dt \Big].  
\enqs

We are then in the framework of Section \ref{Sec:infinitepb} with $d=m=1$ (one-dimensional state variable and control), and the coefficients in \eqref{pb:coeffSDEINF} and \eqref{pb:coeffPayoffINF} are given by (see also Remark \ref{remmultiinfi}):
\begin{equation*}
C = -1, \quad
D^1=\sigma, \quad
N = \delta + \eta, \quad
N + \tilde N = \delta + \eps, \quad
I = - \frac{c}{2}, \quad
H_t =  \frac{cx_0- P_t^0}{2},
\end{equation*}
while the other coefficients are identically zero. Notice that $H_t$ is stochastic.  Under the following assumptions
\begin{equation} \label{condP0}
\int_0^\infty e^{-\rho t} \eee[|P^0_t|^2] dt < \infty, \qquad\qquad \rho > \sigma^2,
\end{equation}
it is clear that {\bf (H1')} and {\bf (H3')} hold true (for the condition in {\bf (H3')}, we can omit the factor $2$, see Remark \ref{RemSimplH3}). The equations for 
$K$ and $\Lambda$ reads as
\begin{equation}\label{eqKLex2}
\left\{
\begin{array}{ccc}
 \frac{ (K+ \frac{c}{2})^2}{\delta+\eta}  + (\rho - \sigma^2) K &=& 0, \\
 \frac{ (\Lambda + \frac{c}{2})^2}{\delta+\eps} + \rho \Lambda - \sigma^2 K &=& 0. 
 \end{array}
 \right.
\end{equation}
Notice that condition {\bf (H2')} is not satisfied. However, we have existence of a solution $(K,\Lambda)$  to  \reff{eqKLex2} such that 
$K_\eta$ $:=$ $\frac{K+c/2}{\delta + \eta}$ $>$ $0$, $\Lambda_\eps$ $:=$ $\frac{\Lambda + c/2}{\delta + \eps}$ $>$ $0$, and given by
\begin{equation}\label{KL}
\left\{
\begin{array}{rcl}
K_\eta  \;: = \; \frac{K+\frac{c}{2}}{\delta+ \eta} &=& \frac{- (\rho-\sigma^2) + \sqrt{ (\rho-\sigma^2)^2 + 2c\frac{\rho-\sigma^2}{\delta+\eta}}}{2}   \; > \; 0, \\
\Lambda_\eps \; := \; \frac{\Lambda + c/2}{\delta + \eps}&=&   \frac{- \rho + \sqrt{ \rho^2 + 2\frac{\rho c + 2\sigma^2 K}{\delta+\eps}}}{2}  \; > \; 0.
\end{array}
\right.
\end{equation}
The equation for $Y$ is written as 
\beqs
dY_t &=& \Big[ \vartheta_t + (\rho + K_\eta)(Y_t  - \E[Y_t]) + (\rho + \Lambda_\eps) \E[Y_t] \Big] dt + Z_t^{2,Y} dW_t^2, \;\;\; t \geq 0, 
\enqs
where $(\vartheta_t)_t$ is the $\F^{W^2}$-adapted process given by
\beqs
\vartheta_t &=&   \frac{K_\eta}{2}(P_t^0 - \E[P_t^0]) + \frac{\Lambda_\eps}{2}(\E[P_t^0] - c x_0). 
\enqs
The solution to the linear mean-field BSDE for $Y$ is explicitly given by 
\begin{equation} \label{expliYex2}
\left\{
\begin{array}{rcl}
Y_t &=&  - \int_t^\infty K_\eta e^{-(\rho+K_\eta)(s-t)} \frac{\E[P_s^0|\Fc_t^{W^2}] - \E[P_s^0]}{2} ds + \E[Y_t], \\
\E[Y_t] &=& - \int_t^\infty  \Lambda_\eps e^{-(\rho+\Lambda_\eps)(s-t)}  \frac{\E[P_s^0] - cx_0}{2}   ds, \;\;\; t \geq 0,  
\end{array}
\right. 
\end{equation}
and clearly satisfies condition {\bf (H4')} from the square integrability condition \reff{condP0} on $P^0$.  
We also notice with Remark \ref{RemSimplH5} that the condition in {\bf (H5')} here writes as $\rho > \sigma^2$, which is satisfied. 
By Theorem \ref{thm:optimalINF},  the optimal control is then given by 
\beqs
\alpha_t^* &=& K_\eta(X_t^* - \E[X_t^*]) + \Lambda_\eps \E[X_t^*] \\
& & \;\; + \; \frac{1}{\delta+ \eta} \Big(Y_t - \E[Y_t] + \frac{P_t^0 - \E[P_t^0]}{2} \Big)  
+  \frac{1}{\delta+ \eps} \Big(  \E[Y_t] + \frac{\E[P_t^0]}{2} - \frac{cx_0}{2} \Big) \\
&=&   K_\eta(X_t^* - \E[X_t^*]) + \Lambda_\eps \E[X_t^*] \\
& & \; + \;  \frac{1}{2(\delta+ \eta)} \Big( (P_t^0 - \E[P_t^0]) -  \int_t^\infty K_\eta e^{-(\rho+K_\eta)(s-t)} (\E[P_s^0|\Fc_t^{W^2}] - \E[P_s^0]) ds  \Big) \\
& & \; + \;   \frac{1}{2(\delta+ \eps)}  \Big( \E[P_t^0] - \int_t^\infty \Lambda_\eps e^{-(\rho+\Lambda_\eps)(s-t)} \E[P_s^0] ds  
- cx_0 \frac{\rho}{\rho+\Lambda_\eps} \Big), 
\enqs
with an optimal level of reserve  given in mean by
\beq
\E[X_t^*] &=& x_0 e^{-\Lambda_\eps t} +   
\frac{\rho c x_0 }{2(\delta+ \eps)}  \frac{1- e^{-\Lambda_\eps t}}{\Lambda_\eps(\rho+ \Lambda_\eps)}  \label{meanXex2} \\
& & \; -  \; \frac{1}{2(\delta+ \eps)}   \int_0^t e^{-\Lambda_\eps(t-s)} \Big(  \E[P_s^0] - \int_s^\infty \Lambda_\eps e^{-(\rho+\Lambda_\eps)(u-s)} \E[P_u^0] du \Big) ds, \;\; t \geq 0. \nonumber
\enq
Suppose that the price $P^0$ admits a stationary level in mean, i.e., $\E[P_t^0]$ converges to some constant $\bar p$ when $t$ goes to infinity (this arises for example when $P^0$ is a martingale or is a mean-reverting process): $\bar p$ is interpreted as a substitute price for the exhaustible resource.  In this case,  it is easy to see from \reff{meanXex2} that the optimal level of reserve also admits a stationary level in mean: 
\beqs \label{explixinfi}
\lim_{t\rightarrow\infty} \E[X_t^*] &=& 
\frac{\rho(cx_0 - \bar p)}{2(\delta+ \eps)\Lambda_\eps(\rho + \Lambda_\eps)} \; =: \; \bar x_\infty.
\enqs
The term $cx_0$ is the cost of extraction for the last unit of resource. When it is larger than the substitute price $\bar p$, i.e.~the Hotelling rent $cx_0 - \bar p$ is positive, 
this ensures that the average long term level of reserve $\bar x_\infty$ is positive, meaning that there is remaining resource when switching to the substitute good. 
This stationary value $\bar x_\infty$ depends a priori on the degree  $\eps$ of  product substitutability, and also on  $\eta$, the penalty parameter for the intermittence of extraction (notice that $\Lambda_\eps$ $=$ $\Lambda_{\eps}(\eta)$ depends on $\eta$ through $K$ $=$ $K(\eta)$, see \reff{KL}). Actually, from straightforward algebraic calculations on 
\reff{eqKLex2}, we have: 
\beqs
2(\delta+ \eps)\Lambda_\eps(\rho + \Lambda_\eps) &=& 
\frac{K_\eta + \rho}{K_\eta + \rho - \sigma^2} (\rho -\sigma^2) c, 
\enqs 
and thus
\beqs
\bar x_\infty &=& \frac{K_\eta + \rho - \sigma^2}{K_\eta + \rho} \frac{\rho}{\rho-\sigma^2} \big( x_0 - \frac{\bar p}{c}\big). 
\enqs
This means that the average long term level of reserve $\bar x_\infty$ $=$ $\bar x_\infty(\eta)$ does not depend on $\eps$, but only on $\eta$. Moreover, as $K_\eta$ is decreasing with $\eta$, this shows that $\bar x_\infty(\eta)$ is also decreasing with $\eta$ (whenever $x_0-\bar p/c$ is nonnegative). 
In particular, $K_\eta$ converges to zero when $\eta$ goes to infinity, and so
\beqs
\bar x_\infty \;  = \; \bar x_\infty(\eta) & \longrightarrow &  x_0  - \frac{\bar p}{c}, 
\;\;\; \mbox{ as } \; \eta \rightarrow \infty. 
\enqs
Finally, notice that the existence of a stationary level of resource in mean implies that  $\lim_{t\rightarrow\infty}\E[\alpha_t^*]$ $=$ $0$: In other words, one stops on average to extract the resource in the long term.


\newpage

\small

\bibliographystyle{plain}

\end{document}